\documentclass[11pt]{amsart}   
\title[Analysis of a free boundary] {Analysis of a free boundary at contact points with Lipschitz data
}

\author[A. L. Karakhanyan]{A. L. Karakhanyan}
\address{Maxwell Institute for Mathematical Sciences and School of Mathematics,
University of Edinburgh, King's Buildings, Mayfield Road, EH9 3JZ,  Edinburgh, Scotland }
\email{aram.karakhanyan@ed.ac.uk}

\author[H. Shahgholian ]{H. Shahgholian\\
 \tiny{\textrm{Accepted in Transactions of the American Mathematical Society}}}
\address{Department of Mathematics, Royal Institute of Technology,  100~44  Stockholm, Sweden}
\email{henriksh@math.kth.se}

\thanks{$2000$ {\it Mathematics Subject Classification.\/} Primary
  35R35. {\it Key words and phrases.\/} Free boundary problem,
regularity, contact points}
\thanks{ H. Shahgholian was partially supported  by the Swedish
Research Council. Authors also thank Professor Carlos Kenig for several valuable comments.
A.K. thanks G\"oran Gustafsson Foundation for visiting appointments to KTH}

\usepackage{graphicx, geometry}

\usepackage{color,fancybox}

\usepackage[dotinlabels]{titletoc}

\usepackage{color}
\usepackage{graphics}

\usepackage{amsmath}
\usepackage{amsfonts}
\usepackage{graphicx}
\usepackage{mathrsfs}
\usepackage{upgreek}
\usepackage{esint} 

\theoremstyle{plain}
\newtheorem{theorem}{Theorem}[section]
\newtheorem{lemma}[theorem]{Lemma}

\newtheorem*{AAA}{Theorem A}
\newtheorem*{BBB}{Theorem B}
\newtheorem*{CCC}{Theorem C}
\newtheorem*{DDD}{Theorem D}
\newtheorem*{EEE}{Theorem E}
\newtheorem{prop}[theorem]{Proposition}
\newtheorem{corollary}[theorem]{Corollary}
\newtheorem{remark}[theorem]{Remark}
\newtheorem{defn}[theorem]{Definition}

\renewcommand\phi\varphi 
\renewcommand\div{\operatorname{div}} 
\newcommand\supp{\operatorname{supp}} 
\newcommand\dist{\operatorname{dist}} 
\newcommand\R{\mathbb{R}} 

\numberwithin{equation}{section}
\newcommand\e{{\varepsilon}}

\newcommand\Rn{{{\mathbb R}^n}}

\newcommand\pr{{\bf Proof:\ }}
\newcommand\p{{\mathcal P}}
\newcommand\la{{\lambda}}
\newcommand\X[1]{{\chi_{\{#1\}}}}

\newcommand\Om{\Omega}
\newcommand{\na}{\nabla}
\newcommand{\outline}[1]{\ifnum\outlineon>0\
\\\noindent\fbox{\begin{minipage}{10cm}{\footnotesize
#1}\end{minipage}}\\\\\fi}
\newcommand{\exclude}[1]{}

\makeindex

\begin{document}

\begin{abstract}
In this paper we consider a minimization problem for the functional
$$
J(u)=\int_{B_1^+}|\nabla u|\sp
2+\lambda_{+}^2\chi_{\{u>0\}}+\lambda_{-}^2\chi_{\{u\leq0\}},
$$
in the upper half ball $B_1^+\subset\R^n, n\geq 2$ subject to  a
Lipschitz continuous Dirichlet data  on $\partial B_1^+$. More
precisely we assume that $0\in \partial \{u>0\}$ and the derivative
of the boundary data has a jump discontinuity. If $0\in
\overline{\partial( \{u>0\} \cap B_1^+)}$ then (for $n=2$ or $n\geq3$ and one-phase case) we prove, among other
things, that the free boundary $\partial \{u>0\}$ approaches the
origin along one of the two possible planes given by
$$
\gamma x_1 = \pm x_2,
$$
where $\gamma$ is an explicit  constant given by the boundary data
and $\lambda_\pm$ the constants seen in the definition of $J(u)$.
Moreover the speed of the approach to $\gamma x_1=x_2$ is uniform.
\end{abstract}
\maketitle


\section{Introduction }

In this paper we consider the local minimizers of the functional
$$
J(u)=\int_{B_1^+}|\nabla u|\sp
2+\lambda_{+}^2\chi_{\{u>0\}}+\lambda_{-}^2\chi_{\{u\leq0\}},
$$
where $B_1^+\subset \R^n, n\geq 2,$ is the upper half of the open unit ball, $\chi_D$ is the
characteristic function of $D\subset \R^n$, $\lambda_\pm$ are given positive constants and $u=f$ on $\partial B_1^+$ with Lipschitz continuous
$f$.
The local regularity of the minimizers $u$  and the free boundary
$\partial \{u>0\} $ were studied in \cite{AC}, \cite{ACF1} and \cite{Gurevich}, notably it was
shown that $u$ is locally Lipschitz continuous.

\smallskip

The boundary regularity of $u$ with  smooth boundary data $f$ such
that $|f(x)|\approx o(|x|)$ near the origin was considered in
\cite{KKS} where, assuming the origin is a contact point, the
authors have proved that close to the origin, the free boundary
approaches the plane $\{x_1=0\}$   in a tangential fashion.

\smallskip

The objective of this paper is to consider boundary data that gives rise to
non-tangential touch between the free and the fixed boundaries.
Such problems appear naturally in the mathematical formulation of the so-called
Dam problem for the water reservoirs (see \cite{AG}). Other problems
of this kind emerge in wake and cavity formations in stationary Eulerian flows moving through
cylindrical domains (see \cite{Birk} Chapters 1.9 and 9.5 for more applications).

Since the formulation of our main results requires some technical
definitions, we refrain ourselves of giving an exact account of our
main results here. However, in lay terms, one can say that our main
result in this paper states that for a boundary data such as $\alpha_+
x_2^+ - \alpha_- x_2^- $, the free boundary $\Gamma (u)$ approaches the fixed one
 along one of the planes $\gamma x_1 =\pm x_2$, where
$$\gamma= \sqrt{\frac{\lambda_+^2-\lambda_-^2}{\alpha_+^2 - \alpha_-^2}-1}.$$
We prove this when $n=2$ for the two phase problem and  $n\geq 3$ for the one phase probelm, see Theorem C. The difficulty for two-phase in higher dimensions
comes from the classification of global homogeneous solutions, that is not feasible by our technique.

\setcounter{tocdepth}{1}
\tableofcontents

\subsection{Plan of the paper}
The plan of this paper is as follows. In this introductory part we
give the necessary notations and definitions to formulate the problem.
Section 2 contains a heuristic discussion of the optimal regularity of
solutions. The key point is the uniform linear growth of minimizers at the
origin. We formulate the main results of this paper in Section 3.
To deal with the boundary behavior of minimizers one needs to obtain up-to boundary
uniform continuity near contact points. The proof of this result as well as a basic
compactness theorem for blow up sequences is contained in Section 4 and Appendix.
Section 5 takes care of the optimal regularity of  minimizers to
our functional.  In Sections 6-8 we show that homogeneous global solutions in one phase case
are two-dimensional, and hence independent of $x_3, x_4, \dots , x_n$.
A stability result is given in Section 9. In fact
Section 9 contains the proof of the main result of this paper,
describing how the free boundary behaves close to the origin. Finally
in Section 10 we give an example of a non-homogeneous global
solution.

\subsection{Notations}
 We will use the following notations throughout the paper.

\begin{tabbing}
$C_0, C_n, \dots$ \hspace{1.29cm}       \=\hbox{generic constants},\\
$\chi_D$            \>\hbox{the characteristic function of the set }$D \subset \mathbb R^n,\  n\geq 2$,\\
$\overline D$       \>\hbox{the closure of } $D$,\\
$\partial D$        \>\hbox{the  boundary of a set }  $D$ ,\\
$x,x'$             \>$x=(x_1,\dots , x_n), \quad x'=(0,x_2, \cdots ,x_n)$,\\
$\R^n_+,\R^n_-$         \>$\{x\in \Rn:\ x_1>0\}; \ \{x\in \Rn:\ x_1<0\}$ , \\
$\Pi$          \>$\{x\in\Rn : \ x_1=0 \}$,\\
$B_r(x), $ $B_r^+(x)$   \qquad    \>$\{y\in \Rn: |y - x|<r\}$, $B_r(x) \cap \R^n_+$ ,\\
$B_r, \  B_r^+$     \>$B_r(0),\  B_r^+(0)$,\\
$B_r^\prime$\>$B_r\cap\Pi$,\\
$S_r^+$\>$\partial B_r\cap\R^n_+$,\\
$\lambda_\pm$, $\Lambda$       \> $\lambda_+, \lambda_-$\ \hbox{ are positive numbers and} {$\Lambda=\lambda_+^2 -\lambda_-^2\neq 0$},\\
$\Gamma(u)$\>         $\partial\{u>0\} $; \hbox{the free boundary of $u$},\\
$\Omega^+(u), \Omega^-(u)$   \>    $\Omega^+(u)=\{x : u(x)>0\},\ \Omega^+(u)=\{x : u(x)<0\}$,\\
$K_\delta(x_0)$   \>   \hbox{the open cone} $K_\delta=\{x\in\R^n_+\ :\ |x-x'|>\delta|x-x_0|\}$,\\
$K_\delta$   \>   \hbox{the open cone} $K_\delta=\{x\in\R^n_+\ :\ x_1>\delta|x'|\}$,\\
${\mathcal P}_r, \mathcal P_\infty, \mathcal{HP_\infty}, \mathcal
P_r'$ \> \hbox{see Definitions
\ref{def-Pr}, \ref{P-classes} and \ref{large-small-stable}},\\
$v^\pm$    \> $v^+=\max(v, 0)$ and $v^-=\max(-v, 0)$. Thus $v=v^+-v^-$.
\end{tabbing}

\subsection{Problem set-up} Throughout this paper we
assume
\begin{equation}\label{conditions-f}
f(x)=\alpha_+x_2^+-\alpha_-x_2^-+g(x),
\end{equation}
where $\alpha_+, \alpha_-$ are nonnegative constants such that  $\alpha_+ + \alpha_- >0$,
and $g(x)\in C^{1,\alpha}(\overline{B_1^+}), \ g(x)=o(|x|)$.
Typically $g(x)=C|x|^{1+\upkappa}$ for positive constants $C$ and $\upkappa$.

\smallskip

For a fixed domain $D\subset\R^n_+$ we put
\begin{equation*}
  J(u, D)=\int_{D}|\nabla u|^2+\lambda_+^2\X{u>0}+\la_{-}^2\X{u\leq
  0}.
\end{equation*}
When it is clear for which $D$ the functional $J$
is considered, we just write it as $J(u)$ omitting the explicit dependence on $D$.
The case $D=B_R^+$ is of particular interest.

\begin{defn}
Let $\mathcal K_f(D)=\{w:w\in H^1(D),\  w-f\in H^1_0(D)\}$
be the class of admissible functions.
\begin{itemize}
 \item A function $u$ is said to be a local minimizer of $J(u, D)$ if for any function $v\in \mathcal K_f(D)$
such that $u=v$ on $\partial D'$, for $D'\subset D$,
it follows that
\begin{eqnarray*}
 J(u)\leq J(v).
\end{eqnarray*}

\item The class of local minimizers is denoted by $\mathcal P(D, n, \lambda_\pm, \alpha_\pm, g)$.
 \end{itemize}
\end{defn}
\begin{remark}\label{rem-p-r}
For  $D=B_r^+$ we denote the corresponding class by
$\mathcal P_r(n,\lambda_\pm, \alpha_\pm, g)$. We also set $\mathcal P_r(n,
\lambda_\pm, \alpha_\pm)=\mathcal P_r(n, \lambda_\pm, \alpha_\pm,
0)$. It is worthwhile to point out that if $u\in \mathcal P_r(n,
\lambda_\pm, \alpha_\pm, g)$ then $u_r(x)=\frac{u(rx)}{r}\in
\mathcal P_1(n, \lambda_\pm, \alpha_\pm, \frac{g(rx )}{r})$ by the
scale invariance of $J(u, B_r^+)$.
\end{remark}

If $D$ is a bounded domain then, from the definition of $J(u, D)$, we have
\begin{equation}\label{J-redef0}
  J(u, D)=\la_{-}^2|D|+\int_{D}|\nabla u|^2+\Lambda\X{u>0},
\end{equation}
where $\Lambda=\lambda_+^2 -\lambda_-^2>0.$ In what follows we take
\begin{equation}\label{J-redef}
 \displaystyle J(u, D)=\int_{D}|\nabla u|^2+\Lambda\X{u>0}.
\end{equation}

Next we introduce  a particular class of local minimizers
$u$, such that the free boundary $\partial\{u>0\}$
is $\delta-${\bf non-tangential} or $\delta-$NT for short.

\begin{defn}\label{def-Pr}
Let $u\in \mathcal P_r(n, \lambda_\pm, \alpha_\pm, g)$.
\begin{itemize}
\item[$\bf 1^\circ$] We say that the free boundary
$\Gamma(u)$ is $\delta-${\textbf{non-tangential}} (or $\delta-{\textbf{NT}}$) at $x_0\in B'_r\cap \Gamma(u)$ if  there exists a $\delta>0$ such
that
\begin{equation}\label{fatness}
  (B^+_{2\rho}(x_0)\setminus B^+_\rho(x_0))\cap\partial\{u>0\}\cap
K_\delta(x_0)\not= \emptyset, \qquad \forall \rho\in(0,r),
\end{equation}
where $K_\delta(x_0)=\{x\in\R^n_+ : x_1>\delta|x'-x_0'|\}$.

\item[$\bf 2^\circ$] The class of all local minimizers in $B_R^+(x_0)$ with $\delta-$NT free boundary is denoted by
$\mathcal P_r(x_0, n, \lambda_\pm, \alpha_\pm, g, \delta)$.
When $x_0=0$ and $R=1$ we often omit the dependence of $\mathcal P_r$ from $x_0$ and write
$\mathcal P_1(n,\lambda_{\pm}, \alpha_{\pm}, g, \delta)$ for brevity.
\end{itemize}
\end{defn}
One can  interpret condition (\ref{fatness}) geometrically as follows:
There is a free boundary point at
each intersection of the cone $K_\delta(x_0)$ with $B^+_{2r}(x_0)\setminus
B^+_r(x_0)$ and hence the free boundary does not approach
the plane $x_1=0$ rapidly as $r\rightarrow 0$.
The next section contains more discussion on
$\delta-$NT as a necessary condition for linear growth.

\smallskip
\begin{remark}
The $\delta-$NT assumption can be weakened as follows. Let
$r>0$ be small, $z\in \partial\{ u>0\}$ be a non-isolated point of the free boundary
and  assume that there is a point $x_r\in (B_{2r}(z)\setminus B_r(z))\cap K_\delta$ such that
\begin{equation}\label{weak-fat}
 |u(x_r)|\leq Cr, \quad \forall r>0
\end{equation}
 for some fixed constants $\delta, C$ independent from $r.$ Then one can prove that
$u$ grows linearly from the origin.
It should be noted here that
(\ref{weak-fat}) is always true for the solutions
to one phase problem provided that the origin is a non-isolated free boundary point, see (\ref{normal-est}).
\end{remark}
\medskip

\subsection{Blow-up limits and Global Solutions}
Let $u_j\in \mathcal P_1(n,\lambda_{\pm}, \alpha_{\pm}, g),\ j=1, 2, \dots$
and  $x_0$ be a contact point, i.e. $x_0\in \Gamma({u_j})\cap B_1'$. Typically $x_0=0$.
For $r_j>0$ we introduce the \textit{blow-up sequence} of functions at $x_0$
\begin{equation} \label{scaling}
v_j(x)=\frac{u_j(x_0+ r_jx)}{r_j}, \qquad r_j\downarrow 0 \ \ {\rm{as}}\ \ j\rightarrow \infty.
\end{equation}
If the sequence $v_j$ is bounded in a suitable space then
sending $r_j$ to 0  we obtain  a so called  \textit{blow-up limit} $u_0$.
One of our main objectives in this paper is to classify the blow-up limits  of the sequence $v_j$ in (\ref{scaling})
as $j$ tends to infinity. It is noteworthy that, in general, the blow-up limit
depends on the sequence $\{r_j\}_1^\infty$. Thus the blow up limit $u_0$
is not unique. Hence it is natural to address the classification of blow up limits.
To do so
we employ the monotonicity formula (\ref{weiss}) and show that the blowup
at the contact points is only one of the functions (\ref{explicit})
(see Sections \ref{Weiss} and  \ref{large-small-sect}).

\smallskip

The classification of all possible blow-up limits is based on geometric properties
that these functions share, notably the linear growth and the homogeneity.

\begin{defn}\label{P-classes} Let $u$ be a local minimizer in $\R^n_+$.
\begin{itemize}
\item[$\bf 1^\circ$] We say that $u$ is a \textbf{global solution} if $u\in\mathcal P_\infty$, where
  \begin{equation*}
  \mathcal P_\infty(C)= \bigcap_{r>0}\left\{u\in {\mathcal P}_r(n, \lambda_\pm, \alpha_\pm)
\ :\ u(0,x')=\alpha_+x_2^+-\alpha_{-}x_2^-, \ |u(x)|\leq
  C(x_1+|x_2|)\right\}
  \end{equation*}
for some positive constant $C$
and $\mathcal P_r(n, \lambda_\pm, \alpha_\pm)=\mathcal P_r(n, \lambda_\pm, \alpha_\pm, 0)$.
\item[$\bf 2^\circ$] The class of all \textbf{homogeneous global solutions} is denoted by
  \begin{equation*}
  \mathcal{HP_\infty}(C)=\left\{u\in\mathcal P_\infty : u(tx)=tu(x),\forall
  t>0\right\} .
  \end{equation*}
\end{itemize}
\end{defn}

This definition requires some explanation. First we note that
any blow-up limit of linearly growing $u$ is a global solution. Moreover
it follows from the monotonicity theorem in Section \ref{Weiss} that
the blow-up $u_0\in\mathcal{HP_\infty}$. The linear growth constant $C$ appearing in the
definition must be consistent with the constants $\alpha_\pm$ that determine the boundary date.
Clearly we must have $C\geq\max(\alpha_+, \alpha_-)$ otherwise at least one of $\alpha_\pm$
must be zero. A posteriori $\mathcal{HP_\infty}$ contains only
two functions, by Theorem C (\ref{explicit}), linking $C$ with constants $\lambda_\pm$ too.
In fact if   $\frac{\lambda^2_+-\lambda_-^2}{\alpha_+^2-\alpha_-^2}-1<0$
then $\mathcal{HP_\infty}$ is empty.  Therefore whenever constant
$C$ is chosen large enough and $\frac{\Lambda}{\alpha_+^2-\alpha_-^2}-1\geq 0$ the
resulted class of global homogeneous solutions is determined uniquely.

\medskip
Finally we define the extreme global solutions and stability
in order to classify the global solutions. 
\begin{defn}\label{large-small-stable}\

\begin{itemize}
\item[$\bf 1^\circ$]$u\in \mathcal P_\infty(n, \lambda_\pm, \alpha_\pm)$ is said to be the smallest (resp. largest)
global solution if for any $v\in \mathcal P_\infty(n, \lambda_\pm, \alpha_\pm)$ we have $u\leq v $ (resp. $u\geq v$).

\item[$\bf 2^\circ$]
The class of all local minimizers that after blow-up
coincide with the smallest homogeneous global solution $v_{S}$
 \begin{equation}\label{P'}
  \mathcal P_r'=\left\{u\in\mathcal P_r(n, \lambda_\pm, \alpha_\pm, g):\lim_{r_j\rightarrow
  0}\frac{u(r_jx)}{r_j}=v_{S}(x), \quad \hbox{for some sequence } r_j\right\}.
  \end{equation}
If $u\in \mathcal P_r'$ then we say that $u$ is \textbf{stable}.
\end{itemize}
\end{defn}

\medskip
\section{Linear Growth: A Heuristic Discussion}\label{Linear Growth}

In analyzing the behavior of the free boundary one needs, in general, to start with the growth rate of the solution at free boundary points.
Lipschitz regularity, up to the boundary, would be the most desirable property for  minimizers of our functional.
This property, or at least the linear growth property at the origin, is indispensable for the rest of the theory to follow.

In general, one cannot expect  this property to hold, and one is forced to impose conditions to assure this.
Indeed, a harmonic function in $B_1^+$ with merely Lipschitz data on $\{x_1=0\}$  is   \textit{not} Lipschitz.
 In such cases  the extra logarithmic term enters into the game, and the solution will belong merely to  the little-o Zygmund class
$$
|u(x)|\leq C |x|\log |x|^{-1}.
$$
In one phase case it is possible to obtain linear growth from the origin, provided  the origin is a non-isolated free boundary point.
In other words if there is a sequence of free boundary points in $\{x_1>0\}$ approaching the origin, then we expect linear growth for the solutions.
 We will state and  give a  proof of this below.
 A similar result of this type was proven in \cite{AG}. Observe that if, even in the one phase case,
we chose the boundary data large enough, e.g. $\alpha_+^2>\Lambda$,  then one may
show that the function $u$ minimizing  $J$ is harmonic in the upper half ball, see Section \ref{large-small-sect}.
Thus, a harmonic function with Lipschitz data can impossibly
be Lipschitz up to the boundary.

\smallskip

For the two phase problem the analysis becomes much more complicated,
and we could not find any complete theory.
Since the Dirichlet data has two signs close to the origin
$$f(x)\approx \alpha_+x_2^+ -\alpha_-x_2^-,$$
the free boundary $\partial\{u>0\}$ is always present in the upper half ball.
The problem is that it might approach the fixed boundary $\{x_1=0\}$ tangentially, and give rise to
a non-Lipschitz behavior of the solution. (This argument does not apply to the one-phase case.)
The  reader may verify that if the free boundary (in two phase case) approaches
tangentially to the fixed boundary and at the same time the solution is Lipschitz then a blow up limit would result in the fact
that one of the  phases vanishes but the boundary data is a two-phase data, and
hence  a contradiction would arise.
This, in particular, suggests that for the two phase problem,
a natural condition to impose is that the free boundary does not touch the fixed one
in a tangential fashion.

It is also not too hard to prove that there are certain Lipschitz boundary data, for which the solution is not Lipschitz and
touches the fixed boundary tangentially. For the proof we would need a classification of homogeneous global solution (as in Theorem C).
Suppose $n=2$, then the proof of Theorem C is more or less elementary in this case (see the proof).
If we accept this result, for the moment, we see that  for $\alpha:=\alpha_+=\alpha_-$,
 and $\Lambda >0$ one may conclude
that the solution cannot be Lipschitz. Otherwise, if this was the case, then
a blow-up of the solution would result in a global solution, with linear growth.
Hence the classification theorem, Theorem C,  would then suggest that the solution is $u=\alpha x_2$,
but then the free boundary condition $|\nabla u^+|^2- |\nabla u^-|^2=\Lambda >0$ fails.

From the representation (\ref{explicit}), we also see that if $\alpha_+^2 - \alpha_-^2>\Lambda $,
then again an up to the boundary Lipschitz continuous solution cannot exist.

The question of finding optimal conditions, that assure linear growth for the minimizers from the origin, is still open.
We have partially answered this question in Theorems A and B, below under mild
conditions.


\medskip

\section{Main Results}
In this section we state the main results of this paper.
To begin our analysis we need the optimal growth estimate
for a local minimizer $u$ near the contact points. More precisely
we have to show that $u\in \mathcal P_1(n , \lambda_\pm, \alpha_\pm, g)$
grows linearly away from $z\in \partial \{u>0\}\cap \Pi$.
Clearly we can assume that $z=0$.

\begin{AAA}\label{linear}
Let $u\in \mathcal P_1(n,\lambda_\pm,\alpha_\pm, g)$ and either of the following holds:
\begin{itemize}
\item[$\bf 1^\circ$]  $u\in \mathcal P_1(n, \lambda_\pm, \alpha_\pm, g, \delta)$,
i.e. the condition (\ref{fatness}) (or its weaker form (\ref{weak-fat}))
is satisfied  for some $\delta >0$ and all $r<1$.

\item[$\bf 2^\circ$] $\alpha_-=0$, $g\geq 0$ and the origin is a non-isolated free boundary point.
\end{itemize}
Then
\begin{equation} \label{linear-growth}
|u(x)|\leq C|x|, \qquad x\in B_{\frac12}^+,
\end{equation}
where $C$ depends on $n,$ $\lambda_\pm,\alpha_\pm,\sup_{B_1} |u|$ and $\delta, g$.
\end{AAA}

\smallskip

As for part $\bf 2^\circ$ of Theorem A, let us note that
the weak $\delta-$NT assumption (\ref{weak-fat}) is always satisfied
for one phase problem, see (\ref{normal-est}).

\smallskip
Our next result is an improvement of Theorem A in the following sense:
Let $u_0$ be a blow-up of $u$ at the origin then
$|u_0(x)|\leq C|x|$ in $\R^n_+$ and $u_0(x)=\alpha_+x_2^+-\alpha_-x_2^-$
on $\Pi$. However these is not enough to conclude that $u_0\in \p_\infty$
since the estimate $|u_0(x)|\leq C(x_1+|x_2|)$ in the definition of
$\p_\infty$ does not follow immediately.
Suppose $T_{i,R}(x)=x+Re_i, i\not =2$ is the translation
in $e_i$ direction by $R\in \R$. Then $u_0(T_{i, R}(x))$
is also a minimizer, but possibly with different constant $C$ in the
linear growth estimate. Does the
boundary data $\alpha_+x_2^+-\alpha_-x_2^-$, depending only on  $x_2$, has any effect?
Do we  get the same growth for $u_0(T_{i, R}(x)  $?

\begin{BBB}\label{Thm B}
Let $u\in \mathcal P_1(n, \lambda_\pm, \alpha_\pm, g)$ and suppose that there is $C>0$
such that
\begin{equation}
|u(x)|\leq C|x-z|,\quad \forall z\in \partial \{u>0\}\cap B_{\frac 12}.
\end{equation}

Then for any blow up limit $u_0$ of $u$ at the origin we have
  $$|u_0(x)|\leq C(x_1+|x_2|).$$
In particular any blow up limit of $u$ belongs to $\mathcal P_\infty(C).$
\end{BBB}

Theorem B is used to classify homogeneous global solutions by employing a
 customary dimension reduction argument. Notably
we show that if $u\in \mathcal{HP_\infty}$ then $u$ depends only on
$x_1$ and $x_2$ variables. Again we note that the growth estimate
$u(x)\leq C |x-z|$ is true for one phase case. As for the two phase
case, one can prove that the uniform $\delta$ or weak
$\delta-$NT condition (see (\ref{weak-fat})) for each contact   point $z\in B_{\frac
12}$ will imply $|u(x)|\leq C|x-z|$ in view of Theorem A.

\smallskip
To set forth the implications of
Theorem B we return to the translated solution $u_0(T_{i,R}(x)), i\geq 3$.
For arbitrary $R_1<R_2$ one can show that  $\max(u_0(T_{i, R_1}(x)), u_0(T_{i, R_2})(x))$ is a
minimizer of $J(u, B_1)$ with boundary values $\max(u_0(T_{i, R_1}(x)), u_0(T_{i, R_2})(x))$
on $\partial B_1^+$. Moreover
by Theorem B the maximum of solutions has exactly the same linear
growth as $u_0$. Thus we can construct a translation invariant maximal global solution.
Repeating this argument for all $i\geq 3$ we obtain a maximal global solution
depending on $x_1$ and $x_2$ only. The minimal solution is
constructed analogously.
Writing Laplace operator in polar coordinates we obtain the classification of
global homogeneous solutions.

\begin{CCC}\label{Thm C}
 In $\R^2$,  there are only two homogeneous global solutions:
  \begin{eqnarray}\label{explicit}
\begin{array}{lll}
v_{L}=\alpha_+(\gamma x_1+x_2)^+  -\alpha_-(\gamma x_1+x_2)^-,\\
  v_{S}=\alpha_+(-\gamma x_1 +x_2)^+ -\alpha_-(-\gamma x_1+x_2)^-,
\end{array}
\end{eqnarray}
where $\displaystyle \gamma =\sqrt{\frac{\Lambda}{\alpha_+^2-\alpha_-^2}-1}.$
Thus $\mathcal{ HP_\infty}=\{v_S, v_L\}$.

This also holds in $\R^n$, for $n > 2$, and for one-phase case, with $\alpha_-=\lambda_-=0$.
 If $\Lambda\leq{\alpha_+^2-\alpha_-^2}$ then there is no free boundary.
\end{CCC}

An obvious consequence of this theorem is that for
any $u\in {\mathcal P}_r$, the angle of the touch between the free
and fixed boundaries is dictated by
the behavior of  $v_S$ or $v_L$.

From Theorem C one can  deduce that the free boundary approaches
the origin along the
 plane $\{x\in \R^n : \gamma x_1= x_2\}$.
The approach is  uniform for the small solution, but in general not for the large
one. For the precise  formulation we introduce some notations: Let
$\sigma$ be a modulus of continuity  and consider

\begin{eqnarray}\label{cones}
&\begin{array}{lll}
  K_\sigma^+ :=\left\{x: \, x_1 >0, x_2 >0, \,
\frac{ x_2}{\gamma+\sigma(|x|)}<x_1<\frac{ x_2}{\gamma-\sigma(|x|)}\right\}, \\
K_\sigma^- :=\left\{x: \, x_1 >0, x_2 < 0, \,
\frac{ -x_2}{\gamma+\sigma(|x|)}<x_1<\frac{- x_2}{\gamma-\sigma(|x|)}\right\},
\end{array}
\end{eqnarray}
\smallskip

\begin{DDD}\label{Thm D}
 Let $u\in  \mathcal P_r$ (see Section \ref{rem-p-r}), and $v_S, v_L$ be defined by (\ref{explicit}).
We consider  $n=2$  for the two phase problem and $n\geq 3$ for the one phase problem.
Then, close to the origin, $\Gamma(u)$ touches  tangentially  one  of the hyperplanes
 $\Gamma({v_S})=\{x\in \R^n: x_2=\gamma x_1\}$ or
$\Gamma({v_L})=\{x\in \R^n: x_2=-\gamma x_1\}$. More precisely
 there exists a modulus of continuity $\sigma(r)=\sigma(u, r)$ and $r_0\in(0,1)$ such that
for any $r\in(0,r_0)$ either
$$\Gamma(u)\subset B_{r}^+\bigcap K_\sigma^+ \qquad or \qquad  \Gamma(u)\subset B_{r}^+\bigcap K_\sigma^-.$$
If $u$ touches the hyperplane $\Gamma({v_S})$ (i.e. $u\in \mathcal P'_1$),
then $\sigma(r)$  and $r_0$ are independent of $u$, and thus the neighborhood $B_{r_0}$ is uniform.


\end{DDD}

\smallskip
It follows from the definition of $\mathcal P_\infty$, and by Theorem B,
that  for $u \in \mathcal P_r$, the limit $u_j(x)=\frac{u(r_jx)}{r_j}, r_j\downarrow 0$
is a global solution.  Furthermore, from Weiss' formula \cite{W}, we have
that the limit has to be a homogeneous function of degree one.
Thus the blow up limits belong to $\mathcal{HP_\infty}$.
However the class of global solutions $\mathcal P_\infty$ may contain non-homogeneous solutions, as our last theorem shows.

\begin{EEE}
 There exists a non-homogeneous global solution with boundary values $\alpha_+x_2^+$.
\end{EEE}

A consequence of Theorem E is a kind of instability of the angle of touch, which amounts to the fact that if a free boundary is asymptotically close to $v_S$, then by
 slight perturbation  of the boundary data the free boundary may come close to $v_L$, asymptotically. This constitutes the idea in the construction
 of global non-homogeneous solutions in Theorem E.

\medskip
Theorem E exhibits the structure of the class of global solutions,
namely the fact that there exist non-homogeneous functions in
$\mathcal P_\infty$. This is due to the following: If $u_j\in
\mathcal P_1(n, \lambda_\pm, \alpha^j_\pm)$ then the blow-up
sequence $v_j=\frac{u_j(r_jx)}{r_j}$ converges to a global solution
$v_\infty\in \mathcal P_\infty(n, \lambda_\pm, \alpha_\pm^\infty)$
where $\alpha_\pm^\infty=\lim_{j\rightarrow \infty}\alpha_\pm^j$.
But it does not necessarily imply that $v_\infty$ is homogeneous. If
$u=u_j$ and $\alpha_\pm=\alpha_\pm^j$ then from Weiss monotonicity
theorem it follows that $v_\infty$ is homogeneous, see Section
\ref{Weiss}.

\bigskip

\section{Technicalities}
In this section we gather a number of useful properties that all
local minimizers share. Some of these properties are of local nature
and some hold true near the fixed boundary, e.g. H\"older continuity.
Although the boundary extensions follow from standard
techniques we have supplied the proofs for the readers' convenience.

\smallskip

\subsection{Uniform H{\"o}lder continuity for $u\in \mathcal P_1(n, \lambda_\pm, \alpha_\pm, g)$}

We begin with recalling some well-known facts, which can be found
can be found in \cite{ACF1}.
\begin{prop}\label{prop}
 Let $u$ be a local minimizer of $J(u)$ in $B_1^+$ and $\Lambda=\lambda_+^2-\lambda_-^2>0$. Then
\begin{itemize}
\item[$\bf 1^\circ$] $u$ is a bounded subharmonic function in $B_1^+$,  Theorem 2.3 \cite{ACF1},
\item[$\bf 2^\circ$] $u$ is harmonic in the interior of $B_1^+\setminus\{u=0\}$,   Theorem 2.4 \cite{ACF1},
\item[$\bf 3^\circ$] $u^+$ is non-degenerate,  Corollary 3.2 \cite{ACF1},
\item[$\bf 4^\circ$] if ${\rm{meas}}\{u=0\}=0$ then $|\nabla u^+|^2-|\nabla u^-|^2=\Lambda$
across the free boundary $\Gamma(u)$ in some weak sense, Theorem 2.4
\cite{ACF1}.
\end{itemize}
\end{prop}

The starting point in our study is the uniform H\"older continuity of local
minimizers.
It will allow us to translate some of the well-known local properties of
$u$ into boundary case.

\begin{lemma}\label{sup-bound}
 Let $u\in \mathcal P_1(n, \lambda_\pm, \alpha_\pm, g)$. Then
$u$ is bounded in $B_{\frac 12}^+$.
\end{lemma}
\pr By Theorem 2.1 \cite{ACF1} $u$ is continuous in each subdomain
$D\subset \subset B_1^+$.  Moreover by Proposition  \ref{prop} $u$ is harmonic
in $\{u\not =0\}$, hence $u^+$ is subharmonic.
Indeed, if $x\in \Om^+(u)$ then $\fint_{B_{r}(x)}u^+\geq u(x)$ for each $r<r_0$ such that
$B_{r_0}(x)\subset\Om^+(u)$,
otherwise  $\fint_{B_r(x)}u^+\geq 0=u^+(x)$ for $x\not \in \Om^+(u)$.
Thus the mean value property is satisfied locally.
Thus $u^+$ is subharmonic.

\smallskip

Let $v$ be the harmonic lifting of $u$, i.e.  $\Delta v=0,
v|_{\partial B_1^+}=u^+$.
From  maximum
principle $u^+\leq v$ and $\int_{B_1^+}|\nabla v|^2\leq
\int_{B_1^+}|\nabla u^+|^2$. In particular
$\|v\|_{H^{1}(B_1^+)}\leq C \|u\|_{H^{1}(B_1^+)}$ with some tame
constant $C$. This yields that $v\in C^{0}(\overline{B_{\frac
12}^+})$. Hence $u^+$ is bounded in $B_{\frac 12}^+$.

By a similar argument one can show that $u^-$ is bounded.
\qed

\medskip

Next theorem is more general and can be applied to families of local minimizers.

\begin{prop}\label{holder-0}
Let $u\in \mathcal P_{R_0}(n, \lambda_\pm, \alpha_\pm, g)$ and
$$\sup_{B_{2R}^+}|u|+\alpha_++\alpha_-+\Lambda+\|f\|_{C^{0,1}}\leq M, \qquad 2R<R_0.$$
Then there are positive constants $\beta=\beta(n, R, M)$ and $C=C(n, R, M)$ such that
$u\in
C^\beta(\overline{B_{R}^+})$ and $\|u\|_{C^{\beta}(B_R^+)}+\|u\|_{H^1(B_R^+)}\leq C$.
\end{prop}

\pr Let $w$ be the harmonic lifting  of $u$ in $B_{2R}^+$.
Because $u-w\in H^1_0(B_{2R}^+)$ then it follows
$$\int_{B_{2R}^+}|\na u|^2-|\na w|^2=\int_{B_{2R}^+}|\na(u-w)|^2+\int_{B_{2R}^+}2\na w\cdot\na(u-w)
=\int_{B_{2R}^+}|\na(u-w)|^2$$
Then from
$J(u, B_{2R}^+)\leq J(w, B_{2R}^+)$ and the equality above we obtain
\begin{eqnarray}\label{M-1}
 \int_{B_{2R}^+}|\na(u-w)|^2&=&\int_{B_{2R}^+}|\na u|^2-|\na w|^2\leq \int_{B_{2R}^+}
\Lambda\left[\X{w>0}-\X{u>0}\right]\\\nonumber
&\leq&
\Lambda|B_1|(2R)^n.
\end{eqnarray}

Take $\eta\in C_0^\infty(B_{2R}), \eta\equiv 1$ in $B_R$, $0\leq \eta\leq 1$
and $|\na \eta|\leq \frac CR$ for some
dimensional constant $C$. Obviously  $(w-f)\eta^2\in H^1_0(B_{2R}^+)$ can be used as
a test function  in the weak formulation of $\Delta w=0$
\begin{eqnarray*}
 \int_{B_{2R}^+}\eta^2|\na w|^2=\int_{B_{2R}^+}\na w\left[\na f\eta^2-2\eta\na \eta(w-f)\right].
\end{eqnarray*}
Applying Cauchy-Schwarz inequality and the estimate $|\na \eta|\leq \frac C R$ we obtain
Caccioppolli's inequality
\begin{eqnarray}\label{M-2}
 \int_{B_R^+}|\na w|^2\leq 8\int_{B_{2R}^+}\left[|\na f|^2+
\frac{4C^2}{R^2}(w-f)^2\right]\leq C_1
\end{eqnarray}
where $C_1=4|B_2|M^2\left[(2R)^{n}+16C^2(2R)^{n-2}\right].$

\smallskip

Since $u-f=0$ in $B_{2R}'$ we can apply Poincar\`e's inequality
to conclude $\int_{B_{R}^+}(u-f)^2\leq\frac{c_0}{R^2}\int_{B_R^+}|\na(u-f)|^2$
depends on the dimension $n$ and $\mathcal H^{n-1}(B_R')-$ the $n-1$ dimensional
Hausdorff measure of $B_R'$.

Combining inequalities (\ref{M-1}), (\ref{M-2}) and Poincar\`e's inequality
we get
\begin{eqnarray}
\int_{B_R^+}|\na u|^2&\leq &2\left(\int_{B_R^+}|\na w|^2+\int_{B_R^+}|\na(w-u)|^2\right) \\\nonumber
&=&2(C_1+\Lambda |B_1|(2R)^n)\equiv C_2
\end{eqnarray}
thereby
\begin{eqnarray}
 \int_{B_R^+}u^2&\leq& 2\int_{B_R^+}f^2+2\int_{B_R^+}(u-f)^2\leq 2\left(M^2\frac{|B_1|}2 R^n
+\frac{c_0}{R^2}\int_{B_R^+}|\na(u-f)|^2\right)\\\nonumber
&\leq & M^2|B_1|R^n+\frac{4c_0}{R^2}\left(\int_{B_R^+}|\na u|^2+\int_{B_R^+}|\na f|^2\right)\\\nonumber
&\leq & M^2|B_1|R^n+\frac{4c_0}{R^2}(C_2+M^2\frac{|B_1|}{2}R^n)\equiv C_3
\end{eqnarray}
implying that $\| u\|_{H^1(B_R^+)}\leq \sqrt{C_2+C_3}\equiv C_4$.

\smallskip

As for H\"older continuity let us note that in view of  Theorem 7.19 of \cite{GT}
it is enough to show that for $B_r^+(z)\subset B_{2R}^+, z\in B_R', r<\frac 12 $
we have
\begin{eqnarray}\label{int-cont}
 \int_{B_r^+(z)}|\nabla u|\leq C_5r^{n-1+\beta},
\end{eqnarray}
for some $\beta>0$ and $C_5$ depending on $M, n$ and $R$.
Indeed if $z\in B_R^+$ and $|z-z'|>\frac 14$ we get that $B_{\frac 18}(z)\in B_1^+$
and by local continuity Theorem 2.1 \cite{ACF1} $u$ is uniformly continuous
with some $\beta>0$ depending only on $\| u\|_{H^1(B_R^+)}, n$ and $M$.
Whilst for $r<\frac12$ either $|z-z'|\leq r$ and
$B_r^+(z)\subset B_{2r}^+(z')$  or $r<|z-z'|<\frac 12.$

\smallskip

First we deal with the case $z\in B_R'$ and $B_{4r}^+(z) \subset
B_R^+$.
Let $v$ be the
harmonic lifting of $u$ in $B_{4r}^+(z)$, i.e. $\Delta v=0$ in $B_{4r}^+$ and  $v-u\in
H^1_0(B_{4r}^+(z)).$ Since $J(u, B_{4r}^+)\leq J(v, B_{4r}^+)$ it follows that
$$\int_{B_{4r}^+(z)}|\nabla u|^2+\Lambda\X{u>0}\leq\int_{B_{4r}^+(z)}|\nabla v|^2+\Lambda\X{v>0}.$$
Thereby
\begin{eqnarray}\label{mean}
\int_{B_{4r}^+(z)}|\nabla u|^2-|\nabla
v|^2&=&\int_{B_{4r}^+}(\nabla u-\nabla v)(\nabla u+\nabla v)\\\nonumber
&=&\int_{B_{4r}^+(z)}|\nabla (u-v)|^2\\\nonumber
&\leq& \int_{B_{4r}^+(z)}\Lambda\X{v>0}-\Lambda\X{u>0}\\\nonumber
&\leq& M|B_4|r^n.
\end{eqnarray}

From triangle inequality we get
\setlength\arraycolsep{2pt}
\begin{eqnarray}\label{grad-integ-est}
\int_{B_r^+(z)}|\nabla u|&\leq& \int_{B_r^+(z)}|\nabla
(u-v)|+\int_{B_r^+(z)}|\nabla v|\\\nonumber
&\leq&
M|B_4|r^n+\int_{B_r^+(z)}|\nabla v|,
\end{eqnarray}
where the last line follows from (\ref{mean}) and Cauchy-Schwarz inequality.

It remains to show that
that there are  constants  $\beta\in(0,1), C_6$ depending on $M, R$ and $n$ such that
\begin{eqnarray}\label{meann}
\int_{B_r^+(z)}|\nabla v|^2\leq C_6r^{n-2+2\beta}.
\end{eqnarray}

To see this  take $\eta\in C_0^\infty(B_{4r})$
such that $\eta\equiv1$ in $B_r$, $0\leq \eta \leq 1, |\na \eta|\leq \frac{C}{r}$,
 $C$ is a dimensional constant,
then
$\eta^2(v-f)=0$ on $\partial B_{4r}^+$ and we have from the weak formulation of harmonicity of $v$
$$\int_{B_{4r}^+}\nabla v[2\eta \nabla \eta  (v-f)+\eta^2 (\nabla v-\nabla f)] =0.$$
Rearranging the terms and applying H\"older inequality we get
\begin{eqnarray*}
 \int_{B_{4r}^+}\eta^2|\nabla v |^2&=&
-\int_{B_{4r}^+}\nabla v\eta[2\nabla \eta  (v-f)-\eta \nabla f] \\\nonumber
&\leq&\e\int_{B_{4r}^+}\eta^2|\nabla v|^2+\frac1\e\int_{B_{4r}^+}[2\nabla \eta  (v-f)-\eta \nabla f]^2.
\end{eqnarray*}
Choosing $\e$ suitably small and recalling that $\eta\equiv 1$ in $B_r$ we get the estimate
\begin{eqnarray}\label{Caccio}
\int_{B_{r}^+}|\nabla v |^2\leq\frac C\e\int_{B_{4r}^+}[2\nabla \eta  (v-f)-\eta \nabla f]^2.
\end{eqnarray}
According to Lemma 1.2.4
in \cite{Kenig} $v$ is H\"older continuous with some exponent $\gamma=\gamma(n, M, R)\in(0,1)$,
because
$|v|\leq M, \|v\|_{H^1(B_{4r}^+)}\leq M+\|u\|_{H^1(B_{4r}^+)}$. Thus the
left hand side of (\ref{Caccio}) can be estimated as follows
$$\int_{B_{4r}^+}[2\nabla \eta  (v-f)-\eta \nabla f]^2\leq
{C_7\sup_{B_{4r}^+}|v-f|}{r^{n-1}}+C_7\sup_{B_{4r}^+}|\nabla
f|r^n\leq C_8r^{n-1+\gamma}$$ where $C_8$ depends only on $n, M, R$
and to get the first inequality we used the estimate $|\nabla
\eta|\leq \frac Cr$. Thus choosing $\beta=\frac{1+\gamma}2$ the
result follows. Notice that $\beta$ depends only on $n, M$ and $R$.

\smallskip

Finally it remains to show  (\ref{int-cont}) for $B_r^+(z)$
with  $z\in B_R^+$ and $r\leq |z-z'|\leq \frac 12.$
Notice that (\ref{grad-integ-est}) and
(\ref{Caccio}) still hold for this case. As for the estimate (\ref{meann}),
it follows from Poisson
representation  and the bound $|v|\leq M$.
\qed

\smallskip
\begin{remark}
 One can apply Proposition \ref{holder-0}  to a countable family
$\mathcal P_{R_j}(n,  \lambda_\pm^j, \alpha_\pm^j, g_j), j=1, 2, ...$ as $R_j\rightarrow \infty$, see the proof of
(\ref{Hol-v_j}) and (\ref{Sob-v_j}) below.
\end{remark}

\subsection{Implications of linear growth }
The standard regularity result for free boundary problems states that the free boundary is
smooth away from an ineluctable singular set of smaller co-dimension.
The genus of regular points is characterized by
flatness.

Mathematically the blow-up consists of scaling
$u$ in small balls centered on the free boundary:
for  $u\in\mathcal P_1(n, \lambda_\pm, \alpha_\pm, g)$ with linear growth at the origin,
the scaled functions $v_j(x)=\frac{u(r_jx)}{r_j}$
are uniformly bounded as $r_j\searrow 0$.
Since $f(0)=0$, one readily verifies that    $v_j\in \mathcal
P_{1/r_j}(n, \alpha_\pm, \lambda_\pm, g_j)$, where $g_j(x)=\frac{g(r_jx)}{r_j}$.
Clearly $v_j$ is defined in $B_{\frac1r_j}^+$ and provides better picture
of the free boundary at the origin.
Thus by scaling we obtain a sequence of function $v_j$ and a sequence of corresponding
free boundaries $\Gamma_j=\Gamma(v_j)$.
One expects that the convergence $v_j\rightarrow v_0$ implies
$\Gamma_j\rightarrow \Gamma_0=\Gamma(v_0)$ in Hausdorff
distance, which will follow immediately from a compactness of $v_j$ in a suitable class of functions.
For the reader's convenience we recall Theorem 3.1 from \cite{KKS}.

\begin{prop}\label{technical}(\cite{KKS})
Let $v_j$ be a blow up sequence of $u_j$, as in (\ref{scaling}),
  with $u_j\in \mathcal P_{1}(n, \lambda_\pm,\alpha_\pm, g) $ and $x_0=0$.
Further assume that $u_j$ have uniform linear growth. Then, after passing to a subsequence, there exists
 $v\in\mathcal P_{\infty}$
so that
\begin{itemize}
\item[$\bf{1^\circ}$]$v_j\rightarrow v$ uniformly on compact subsets of $\R^n_+$ and
        in $C^{\beta}(E), 0<\beta<1,$ for each $E\subset\subset\R^n_+$,
\item[$\bf{2^\circ}$]for each $M$, $v_j\rightharpoonup v$ weakly in $H^{1}(B_M^+)$,
\item[$\bf{3^\circ}$]for each $M$, $\chi\{v_j>0\}\rightarrow\chi\{v>0\}$
in $L^1(B_M^+)$,
\item[$\bf{4^\circ}$]$\nabla v_j(x)\rightarrow\nabla v(x)$ for a.e. $x$,
\item [$\bf{5^\circ}$] For each $\delta >0$, $E\subset B^+_M$, $dist (E,\Pi) \geq
  \delta$, $0<r<\delta/4$, for $j$ large
$$\partial \{v_j>0\} \cap E \subset \bigcup\limits_{x\in \{v>0\}\cap
  E_{\delta/2}} B_r(x),
$$
and
$$
\partial \{v>0\} \cap E \subset  \bigcup\limits_{x\in \{v_j>0\}\cap
  E_{\delta/2}} B_r(x),
$$
where $E_{\delta/2}$ is a $\delta/2$-neighborhood of $E$.
\end{itemize}
\end{prop}

\subsection{Weiss' energy}\label{Weiss}
It follows from \cite{W} (see also \cite{W-bndr})
that  for  any $u\in \mathcal P_r(n, \lambda_\pm, \alpha_\pm, 0)$

\begin{equation}\label{weiss}
  W(R, u, x_0)=W(R)=\frac1{R^n}\int_{B_R^+(x_0)}|\nabla u|^2+\Lambda \X{u>0}-
\frac1{R^{n+1}}\int_{S^+_R(x_0)}u^2,
\end{equation}
is non-decreasing function of $R$, with $x_0\in \Gamma(u), B_R(x_0)\subset B_r^+$,
and
$$\frac{dW}{dR}=\frac1{R^{n}}\int_{S_R^+}\left(\nabla u\cdot \nu-\frac uR\right)^2.$$
 $W(R, u, x_0)$ is called Weiss' energy at $x_0$.
Notice that $\nabla u\cdot x-u=0$ if and only if $u$ is  homogeneous functions of degree 1.

\begin{prop}\label{W-prop}
Let $u\in \mathcal P_1(n, \lambda_\pm, \alpha_\pm,  g)$ and
$g(x)=C|x|^{1+\upkappa}$. If $u$ has linear growth
then $W(R, u, 0)$ is non-decreasing function of $R$
and
$$\frac{dW}{dR}\geq \frac1{R^{n}}\int_{ S_R^+}\left(\nabla u\cdot \nu-\frac uR\right)^2.$$
In particular any blow-up limit of $u$ at the origin is homogeneous function of degree one.
\end{prop}

\begin{remark}
 To the benefit of clarity we take $g(x)=C|x|^{1+\upkappa}$ with $C, \upkappa>0$.
The case of more general $g(x)=o(|x|)$ can be dealt with similarly, namely
one needs to add a corrective term to $W$ to maintain the monotonicity.
\end{remark}

\pr
If $g\not =0$ and $u\in \mathcal P_r(n, \alpha_\pm, \lambda_\pm, g)$
then some extra care is needed to prove the estimate from below for the derivative $W'(R, u, 0)$.
See Lemma \ref{mon-aux} in Appendix for the proof.

\medskip

It remains to show that $W(r, u, 0)$ is bounded when
$r$ tends to zero. If  $v$ is the harmonic lifting of $u$
in $B_{4r}^+$ and  $u$ has linear growth at
$0$, i.e. $\sup_{B_{4r}^+}|u|\leq  Cr$, then by maximum principle  $\sup_{B_{4r}^+}|v|\leq Cr$.
From Caccioppolli's inequality (\ref{Caccio}) we have

$$ \int_{B_r^+}|\nabla v|^2\leq  Cr^n.$$
Hence
$$\int_{B_r^+}|\nabla u|^2\leq 2\int_{B_r^+}|\nabla v|^2+2\int_{B_r^+}|\nabla (u-v)|^2$$
which, in view of  (\ref{mean}), implies
that $W$ is bounded for small $r$, whenever
$u\in \mathcal P_1(n, \lambda_\pm,  \alpha_\pm, g)$ is linearly growing solution.
\qed


\bigskip
\section{Proof of Theorem A}
The proof of Theorem A consists
of two parts. The first one deals with the two phase problem. Our
method is based on dyadic scaling argument. If the statement of Theorem A
fails then it allows us to construct a linearly growing, non-degenerate harmonic function $v_0$
in $\R_+^n$ vanishing on $\partial \R^n_+$ and at some interior point of $\R^n_+$.
The latter is
due to $\delta-$NT condition, see \ref{fatness}. Thus, in view of the Liouville theorem,
$v_0$ is zero, which contradicts the non-degeneracy of $v_0$.

\subsection{Two-phase case}
Set
$$
S(j,u):=\sup_{B_{2^{-j}}^+}|u|.
$$
It suffices to show
\begin{eqnarray}\label{discrete-linear}
 S(j+1,u)\leq\max\left\{\frac{c2^{-j}}{2},\frac{S(j,u)}{2},\dots,\frac{S(0,u)}{2^{j+1}}\right\}
\end{eqnarray}
for some positive constant $c$. Let us suppose that (\ref{discrete-linear}) is not true.
Then there exists %
a sequence of minimizers $u_j\in\mathcal P_1(n, \lambda_\pm, \alpha_\pm, g)$  and a sequence of
integers $k_j$ so that
\begin{equation}\label{linear2}
S(k_j +
1,u_j)>\max\left\{\frac{j2^{-k_j}}{2},\frac{S(k_j,u_j)}{2},\dots,
\frac{S(k_j-m,u_j)}{2^{m+1}},\dots,\frac{S(0,u_j)}{2^{k_j+1}}\right\}.
\end{equation}
 Observe that from Lemma \ref{sup-bound}
$|u_j|\leq M$ hence $k_j\rightarrow \infty$. Put
$$v_j(x)=\frac{u_j(2^{-k_j}x)}{S(k_j + 1,u_j)}.$$

We wish to show that (\ref{linear2}) implies uniform up-to-boundary estimates for the sequence
$v_j$. In fact there are positive constants $\alpha$ and $C$
depending on $R$ but independent of $j $ such that the following estimates hold

\begin{eqnarray}
 \label{Hol-v_j} \|v_j\|_{C^\alpha(B_R^+)}\leq C(R),\\\
\label{Sob-v_j} \|v_j\|_{H^1(B_R^+)}\leq C(R).
\end{eqnarray}

For brevity we denote

\begin{eqnarray}\label{j-stuff}
  \epsilon_j=\frac{2^{-k_j}}{S(k_j+1,u_j)},
\qquad f_j(x)=\epsilon_j(\alpha_+x_2^+-\alpha_-x_2^-)
+g_j(x),
\end{eqnarray}
where $g_j(x)=\frac{g(2^{-k_j}x)}{S(k_j+1,u_j)}.$
Recall that by (\ref{linear2})
\begin{eqnarray}\label{sigma-zero}
             \displaystyle\epsilon_j=\frac{2^{-k_j}}{S(k_j+1,u_j)}\leq\frac{1}{j}\longrightarrow 0
            \end{eqnarray}
thereby
$f_j\rightarrow 0$ when $j\rightarrow \infty$, for $g(x)=o(|x|)$.

Consider the
{\em{scaled}} functional \setlength\arraycolsep{2pt}
\begin{eqnarray}\label{scaled-J}
\widetilde J_j(v,B_R^+)=\int_{B_R^+} |\nabla v|^2+\epsilon_j^2\Lambda\X{v>0}.
\end{eqnarray}

If  $u_j\in\mathcal P_1(n,\la_{\pm},\alpha_{\pm},g)$ then
$v_j\in \mathcal P_{2^{k_j}}(n, \epsilon_j\alpha_\pm, \epsilon_j\lambda_\pm, g_j)$
provided $R<2^{k_j}$. Indeed by a simple calculation we have
\begin{eqnarray*}
\widetilde J_j(v_j, B_R^+)&=&\int_{B_R^+}|\nabla
v_j|^2+\epsilon_j^2\Lambda\X{v_j>0}\\\nonumber
&=&\epsilon_j^22^{k_jn}\int_{B_{R/2^{k_j}}^+}|\nabla
u_j|^2+\Lambda\X{u_j>0}\\\nonumber
&=&\epsilon_j^22^{k_jn}J(u_j,B^+_{R/2^{k_j}})\nonumber.
\end{eqnarray*}

Furthermore for fixed $R=2^m$ we infer from (\ref{linear2}) that
  \begin{itemize}
    \item $\sup\limits_{B_{\frac{1}{2}}^+}|v_j|=1$,
    \item $\sup\limits_{B_{2^m}^+}|v_j|\leq C2^m , R=2^m<2^k_j$, $m$ is fixed.
  \end{itemize}

Now we can apply Proposition  \ref{holder-0} with
$\sup\limits_{2R}|v_j|+\epsilon_j(\alpha_++\alpha_-+\lambda_++\lambda_-)+\|f_j\|_{C^{0,1}}\leq M$
with $M=2^{m+1}$ and the estimates
(\ref{Hol-v_j}) and (\ref{Sob-v_j}) follow.

Thereby we can extract a subsequence $v_{j_k}$ which converges to
some function $v_0$ such that the following holds: for any fixed $R>0$

\begin{eqnarray}\label{limit-function}
\left\{
\begin{array}{llllll}
(i) & v_{j_k}\rightarrow v_0\  \mbox{in}\
C^{\beta}(\overline{B^+_R}), \qquad   v_{j_k}\rightharpoonup  v_0\  \mbox{weakly in}\
H^{1}_{loc}(\R^n_+) ,\\
(ii) & \sup\limits_{B_{1/2}^+}|v_0|=1, \  v_0(x)=0, x\in \Pi\qquad \mbox{by}\ C^{\beta}\ \mbox{regularity},\\
(iii)& \Delta v_0=0\  \mbox{in}\  x_1>0, \\
(iv) & v_0\mbox{ has linear growth},\\
(v) & v_0(y_0)=0\  \mbox{for some interior point} \ y_0\   \mbox{(by (\ref{fatness}))}.\\
\end{array}
\right.
\end{eqnarray}

Once all  claims in (\ref{limit-function}) are proven we may use
Liouville's theorem for harmonic functions in $\R^n_+$
(utilizing (iii) and (iv)) to conclude $v_0(x)=ax_1$ for some constant $a\not=0$.
But then (ii), (v) and (vi) are in direct contradiction,
and hence our supposition (\ref{linear2}) is false.

\smallskip

Now we proceed by proving (\ref{limit-function}).
The first  claim follows from  standard compactness arguments. The second
one follows from (\ref{j-stuff}) and the convergence of the
traces of $v_j$ in  view of H\"older continuity.

Let us prove the third claim. Let $D\subset \overline{B_R^+}$ be a
domain and $R>0$ is fixed. Then
$v_j\in \mathcal P_{2^{k_j}}(n, \epsilon_j\alpha_\pm, \epsilon_j\lambda_\pm, g_j)$
for the  {\em scaled} functional
$\widetilde J_j$, defined by (\ref{scaled-J}).
Observe that for each
$\psi\in C^\infty_0(D)$
$$\widetilde J_j(\psi, D)\rightarrow \int_D|\nabla \psi|^2 \qquad\mbox{as}\qquad j\rightarrow\infty .
$$

By  (\ref{Hol-v_j}) and (\ref{Sob-v_j}),  $v_0$ exists and
$\int_{D}|\na v_0|^2\leq\liminf\limits_{k\rightarrow 0}\int_D|\na v_{j_k}|^2 $.
According to  (\ref{j-stuff}) $f_0=v_0=0$ on $\Pi=\{x: x_1=0\}$,
where $f_0=\displaystyle\lim_{j\rightarrow \infty} f_j$ uniformly.

Now let us  take $\psi\in H^1_0(D)$, then
\begin{eqnarray}\nonumber
\widetilde J_j(v_j, D)\leq \widetilde J_j(v_j+\psi, D)
\end{eqnarray}
or equivalently
\begin{eqnarray}\nonumber
 \int_D\epsilon_j^2\Lambda\X{v_j>0}\leq \int_D -2\nabla v_j\nabla\psi+|\nabla \psi|^2+
\epsilon_j^2\Lambda\X{v_j+\psi>0}.
\end{eqnarray}
Thereby sending $j_k$ to $\infty$ and utilizing the weak convergence of gradients
$\nabla  v_{j_k}\rightharpoonup \nabla v_0$ in $L^2(\overline{B_R^+})$,  we conclude
\begin{eqnarray}\nonumber
0\leq \int_D -2\nabla v_0\nabla\psi+|\nabla \psi|^2
\end{eqnarray}
and upon adding $\int_D|\nabla v_0|^2$ to both sides we infer
\begin{eqnarray}\nonumber
\int_D|\nabla v_0|^2\leq \int_D |\nabla(v_0- \psi)|^2.
\end{eqnarray}
Since $C_0^\infty(D)$ is dense in $H^{1}_0(D)$ we conclude
the proof of the third claim in (\ref{limit-function}).

The fourth claim follows from (\ref{linear2}) as indicated above.
Hence it remains to
prove the fifth claim. By our assumption (\ref{fatness}) (resp.  (\ref{weak-fat}))
 there exists $x_j\in B_{2^{-k_j}}^+\cap K_\delta$ such that $u_j(x_j)=0$ (resp. $|u(x_j)|\leq C|x_j|$).
Thereby
$$\frac{1}{2}\leq\frac{|x_j|}{2^{k_j}}\leq 1.$$
If we set $y_j=\frac{x_j}{2^{k_j}}$, then one can easily verified that
$y_j\in (B_1\setminus B_{1/2})\cap K_\delta$ and $y_j\rightarrow y_0$,
for some $ y_0\in (B_1^+\setminus B_{1/2}^+)\cap K_\delta$.
Clearly  $v_0(y_0)=0$ by (\ref{j-stuff}) and H\"older continuity.

Now the proof of (\ref{linear-growth}) for two phase case  is complete.

\subsection{One-phase case}
To prove (\ref{linear-growth}) in Theorem A $\bf 2^\circ$,
we need to work out the condition (v) in  (\ref{limit-function}), because the
others follow as above.
For two-phase case, (v) was justified by  assumption (\ref{fatness}) whilst
for one-phase case, (\ref{fatness}) is replaced by the
condition that the origin is a non-isolated free boundary point.
Indeed, this would be enough to force through a similar condition
as that in (v) of (\ref{limit-function}).
However, the analysis  is slightly more delicate and needs care.

\begin{figure}
\begin{center}
\includegraphics[scale=0.9]{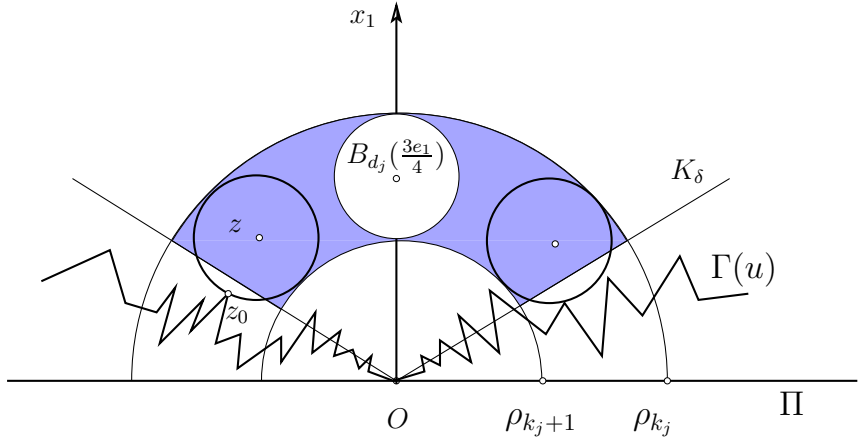}
\caption{The weak $\delta$-NT condition for one-phase problem.}
\end{center}
\end{figure}

\smallskip

Suppose for a sequence  $k_j\uparrow\infty$ we have
$\{u=0\}\cap(B_{\rho_{k_j}}\setminus B_{\rho_{k_j+1}})\cap K_\delta \not = \emptyset$,
where $\rho_{k_j}=\frac 1{2^{k_j}}$.
We consider the family of balls  $B_{d_j}(\zeta)$ with $d_j=\frac14\rho_{k_{j}}$, such that
$\zeta \in S^+_{\frac34\rho_{k_j}}$ and $B_{d_j}(\zeta)$ is above the free boundary
$\Gamma(u)$. Then in this family of balls there is one that touches the free boundary $\Gamma(u)$.
Let $B_{d_j}(z)$ be such ball touching the free boundary
at $z_0$. Clearly $u$ is positive and harmonic inside $B_{d_j}(z)$
and attains its minimum at $z_0$,
therefore we can apply   Lemma 11.19 from \cite{CS} to get  the estimate
\begin{equation}\label{normal-est}
 u(z)\leq C d_j \frac{\partial u}{\partial \nu}(z_0),
\end{equation}
where $\nu$ is the inner normal to $B_{d_j}(z)$ at $z_0$.
Then by  Theorem 6.3 in \cite{AC} we have  $|\nabla u(z_0)|\leq \lambda_+$,
which in conjunction with    Harnack's inequality
implies
$\displaystyle\sup_{B_{\frac{d_j}2}(z)} u\leq C_0 u(z)\leq C_0 C \lambda_+ d_j=\frac{C_0 C\lambda_+}{4}
\rho_{k_j}$. Hence
\begin{eqnarray}\label{0-case-2}
 \sup_{B_{\frac{d_j}2}(z)} u\leq \frac{C_0 C\lambda_+}{4} \rho_{k_j} .
\end{eqnarray}

For scaled functions $v_j(x)=\frac{u(\rho_{k_j}x)}{S(k_j+1, u_j)}$ it follows from (\ref{0-case-2}),
that there exists a ball $B_{\frac 14}(y_0)\subset B_1^+\setminus B_{1/2}^+$ such that
$$\sup_{B_{\frac14}(y_0)} v_j\leq C\frac{\rho_{k_j}}{S(k_j+1, u_j)}=C\epsilon_j\rightarrow 0$$
by (\ref{sigma-zero}) which gives (v) in (\ref{limit-function}) for one phase case.

The proofs of the remaining claims of (\ref{limit-function}) are the same as for the two-phase case
and one will have the final contradictory conclusion.\qed


\bigskip
\section{Proof of Theorem B}

It follows from the proof of Theorem A $\bf 2^\circ$, that $u\geq 0$ grows linearly
away from the origin, provided the origin  is a non-isolated free boundary point.
We can replace the origin by any non-isolated free boundary point
$z$ near the origin and apply the same argument to show that
for $u\in\mathcal P_1(n, \lambda_\pm, \alpha_\pm, g)$
there exists a tame constant $C$ such that the growth estimate
\begin{equation}\label{linear-z}
 0\leq |u(x)|\leq C|x-z|
\end{equation}
holds for any $z\in B_r'\cap\partial\Om^+(u)$ for some $r>0$.

\smallskip

In order to conclude (\ref{linear-z}) for the two
phase solutions we further require the  $\delta-$NT condition to be
satisfied in some neighborhood of the origin. Notice that in the two
phase case, by the H\"older continuity of $u$, the origin  is
automatically a non-isolated free boundary point.

\smallskip

Our goal is to prove that the free boundary $\Gamma(u)$ remains
within a cone $\mathcal C_{\delta_0} =\{x: x_1\geq{\delta_0}|x_2|\}$
in some neighborhood of the origin. This will be enough to prove
Theorem B, because for the free boundary of the blow-up it implies
$\Gamma(u_0)\subset \mathcal C_{\delta_0}.$ Thus the uniform
$\delta-$NT condition will be satisfied for $u_0$, with
${\delta_0}=\delta$ and the result will follow from Theorem A via a
standard scaling argument.

\begin{lemma} Let $u\in \mathcal P_1(n, \lambda_\pm, \alpha_\pm,
g)$. If the $\delta-$NT  assumption (\ref{fatness}) is satisfied for
any free boundary point $x_0\in B_{\frac12}'$ then there exists  a
tame constant $\delta_0$  such that
  \begin{equation}\label{6.1}
  \frac{|z_2|}{z_1}\leq \delta_0, \qquad \forall z\in \Gamma(u)\cap B_{1/2}^+.
  \end{equation}
In particular for any blow-up limit $u_0$ the inclusion
$\Gamma(u_0)\subset \mathcal C_{\delta_0}$ is true.
\end{lemma}

\pr
It follows from the uniform $\delta-$NT condition
and the discussion above that (\ref{linear-z}) is true.
Suppose (\ref{6.1}) fails, then there exists a sequence
$z^k\in B_{1/2}^+\cap \Gamma(u_k)$ of free boundary points of
$u_k\in \mathcal P_1(n, \lambda_\pm, \alpha_\pm, g)$,  such that
$$|z^k_2|\geq kz^k_1,\ k>k_0,$$
for sufficiently large $k_0\in \mathbb N$. 
Setting $d_k=z^k_1, f_0(x)=\alpha_+x_2^+-\alpha_-x_2^-$ we have
\begin{eqnarray*}
|z^k_2|\geq kd_k, \qquad |f_0(p^k)|=\left\{\begin{array}{ll}\alpha_+|z_2^k| \ {\rm{if }}\  z_2^k>0 \\
\alpha_-|z_2^k| \ {\rm{if }}\  z_2^k<0\end{array}\right. ,
\end{eqnarray*}
where $p^k$ is the projection of $z^k$ onto $\Pi$.
In any case we get that $|f(p^k)|\geq k \min(\alpha_+, \alpha_-)d_k.$
Put $r_k=|z^k-\xi^k|$, where $\xi^k=(0, 0, z^k_3,\dots, z_n^k)$ is the projection of $z^k$ onto
$\Pi \cap \{x_2=0\}$. Then from triangle inequality
$ r_k=|z^k-\xi^k|\geq |z^k_2|-z_1^k>(k-1)d_k$ implying $\frac{d_k}{r_k}\leq \frac1{k-1}$. In particular
\begin{equation}\label{r-k}
 1\geq |y^k_2|=\displaystyle\frac{|z^k_2|}{r_k}\geq 1-\frac{d_k}{r_k}\geq1- \frac{1}{k-1} .
\end{equation}

Now introduce the  scaled functions $$v_k(y)=\frac{u_k(\xi ^k+{r_ky})}{r_k}, \qquad y\in B_{2}^+.$$
The points $y^k=\displaystyle\frac{z^k-\xi^k}{r_k}$ are on
the half sphere $S_{1}^+$
and  from (\ref{r-k}) we get
$$y^k_1=\frac{d_k}{r_k}<\frac{1}{k-1}.$$

\begin{figure}
\begin{center}
\includegraphics[scale=0.7]{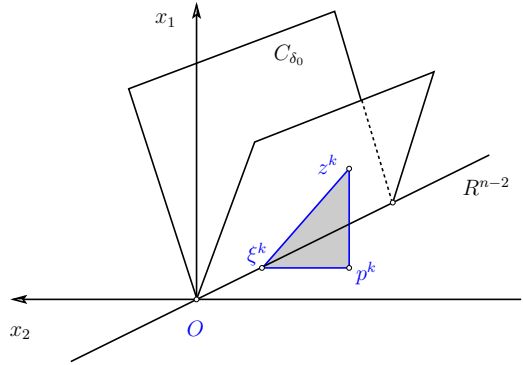}
\caption{The cone $\mathcal C_{\delta_0}$.}
\end{center}
\end{figure}

From (\ref{linear-z}) we have
$|u_k(x)|\leq C|x-\xi ^k|$, since  $\xi^k\in B'_{\frac12}\cap\Gamma(u_k)$.
Therefore it follows that $|v_k(y)|\leq C|y|$, with constant $C$ independent of $k$.
Furthermore $v_k$ is a local minimizer of $J(\cdot, B_2^+)$ because
\setlength\arraycolsep{2pt}
\begin{eqnarray*}
  \int_{B_1^+}|\nabla v_k|^2+{\Lambda}\chi_{\{v_k>0\}}=
\frac{1}{r_k^n}\int_{B_{r_k}}|\nabla
  u|^2+{\Lambda}\chi_{\{u>0\}}. 
\end{eqnarray*}
Thus $v_k\in \mathcal P_2(n,\lambda_\pm, \alpha_\pm, g_k)$
where $g_k(x)=\frac{g(\xi ^k+{r_ky})}{r_k}$ and $\lim g_k=0$ uniformly.

By Proposition \ref{holder-0} it follows,  that
$v_{k}$ is bounded in $C^\beta(\overline{B^+_{3/2}})\cap H^{1}(\overline{B_{3/2}^+})$ for
some positive $\beta\in(0,1)$. Then for a subsequence $k_j$,  $v_{k_j}\rightarrow v_0$ in
$C^\beta(\overline{B^+_{3/2}})$, $\nabla v_{k_j}\rightharpoonup \nabla v_0 $
weakly in $L^2(B^+_{3/2})$
and $y^{k_j}\rightarrow y^0$, where $y^0$ is a free boundary point.
From (\ref{r-k})  $|y^k_2|=\displaystyle\frac{|z^k_2|}{r_k}\rightarrow 1$
and $y\in S_1^+$.
But then
$y^0_1=0, |y_2^0|=1$ and this contradicts to $f_0(y^0)=\alpha_\pm\not=0$.\qed

\smallskip

Let $v(x)=\frac{u_0(\xi+Rx)}{R}$ then  (\ref{6.1}) translates to the
free boundary of the blow-up function $u_0$ implying that
$\Gamma(u_0)\subset \mathcal C_{\delta_0}$. Hence we have uniform
$\delta-$NT condition for each $z\in \Gamma(u_0)\cap \Pi$.
From  Theorem A we have $|v(x)|\leq C|x|$. Returning to $u$ we
conclude $|u(x)|\leq C|z-\xi|\leq C(x_1+|x_2|)$ and this finishes
the proof of Theorem B.

\section{Largest and Smallest global solutions
}\label{build}
Before embarking  into the details we briefly go over the main steps of the proof.
First  we notice that the global solutions enjoy ordering.
This implies that there are  smallest and largest global homogeneous
solutions which we denote  respectively by $v_S$ and $v_L$.
 It follows from the scale and translation invariance that $v_S$ and $v_L$
depend only on $x_1$ and $x_2$. Hence we can explicitly compute them. Moreover
$v_S$ has larger $W$-energy implying that the free boundary of any global homogeneous
solution, distinct from $v_L$ and $v_S$,  cannot touch
$\Gamma({v_S})$ or $\Gamma({v_L})$ tangentially.

\smallskip

Thus if there is third global homogeneous solution $u$ then
we can construct a new one which is symmetric in $x_3, x_4, \dots, x_n$
variables and neither of the functions $v_S, v_L$ coincides with $u$.
Thus  without loss of generality we may assume that $u$ is symmetric in $x_3, x_4, \dots, x_n$
variables.
Then a dimension reduction argument will finish the proof since the only 2D solutions
are $v_S$ and $v_L$.


\subsection{Largest and smallest solutions in $\mathcal P_\infty$}\label{large-small-sect}
We recall (\ref{J-redef})
$$J(u, B_R^+)=\int_{B_R^+}|\nabla u|^2+\Lambda\X{u>0}.$$
Let $v_1, v_2$ be two
minimizers of $J(u, B_R^+)$ and $v_1\leq v_2$ (resp. $v_1\geq v_2$) on $\partial B_R^+$.
Then it is easy to see that $\max(v_1, v_2)$ (resp. $(\min(v_1, v_2)$) is a minimizer of $J(u, B_R^+)$
with boundary values $v_2$ (resp. $v_1$).

Indeed testing  $\max(v_1,v_2)$  against $v_2$ in $B_R^+$ and $\min(v_1,v_2)$
against $v_1$ we get
\begin{eqnarray}\label{star}
  J(v_2, B_R^+)\leq J(\max(v_1,v_2), B_R^+),&\\\nonumber
   J(v_1, B_R^+)\leq J(\min(v_1,v_2), B_R^+).&\nonumber
\end{eqnarray}

Clearly
\begin{eqnarray}\nonumber
J(\max(v_1,v_2), B_R^+)=\int_{B_R^+\cap\{v_1>v_2\}}|\nabla v_1|^2+\Lambda\chi\{v_1>0\}+
\int_{B_R^+\cap\{v_1\leq v_2\}}|\nabla
v_2|^2+\Lambda\chi\{v_2>0\}\\\nonumber
J(\min(v_1,v_2), B_R^+)=\int_{B_R^+\cap\{v_1>v_2\}}|\nabla v_2|^2+\Lambda\chi\{v_2>0\}+
\int_{B_R^+\cap\{v_1\leq v_2\}}|\nabla
v_1|^2+\Lambda\chi\{v_1>0\}\\\nonumber
\end{eqnarray}
which  gives
\begin{eqnarray}\label{2star}
J(\max(v_1,v_2), B_R^+)+J(\min(v_1,v_2), B_R^+)=J(v_1, B_R^+)+J(v_2, B_R^+).
\end{eqnarray}

Hence (\ref{star}) in conjunction with (\ref{2star}) implies
\begin{eqnarray*}
  J(v_1, B_R^+)=J(\min(v_1,v_2), B_R^+),\\\nonumber
  J(v_2, B_R^+)=J(\max(v_1,v_2), B_R^+).
\end{eqnarray*}

Upon applying this observation to finite number of minimizers we obtain

\begin{lemma}\label{rem-max-min}
If $v_1, \dots, v_N$ are
minimizers on $B_R^+$ and $v_1\leq v_2\leq\dots \leq v_N$  on $\partial B_R^+$
(resp. $v_1\geq v_2\geq \dots\geq v_N$)  then $v_{L}^R=\max(v_1, \dots, v_N)$
(resp. $v_{S}^R=\min(v_1, \dots, v_N)$) is a minimizer  of
$J^R$ with  boundary values $v_N$ on $\partial B_R^+$.
\end{lemma}
Employing a  compactness argument it follows that there exists
a largest and a smallest minimizer denoted respectively by $v_L^R$ and $v_S^R$.

By definition,  for
any $u\in\p_R(n, \lambda_\pm, \alpha_\pm)\cap\mathcal P_\infty$ we have
$$
v_{S}^R(x)\leq u(x)\leq v_{L}^R(x),\qquad x\in B_R^+.
$$
Moreover by Definition \ref{P-classes}, $v_S^R$ and $v_L^R$ have uniform linear
growth, i.e. $|v_S^R|, |v_L^R|\leq C(x_1+|x_2|)$ for some tame constant $C$ independent of
$R$. Sending $R\rightarrow\infty$ and  utilizing the linear
growth   Proposition \ref{holder-0}  we infer that $v_{L}^R\rightarrow v_{L}$ uniformly and
weakly in $H^1_{\textrm{loc}}$. Furthermore
$v_{L}\in\p_\infty$. \

Indeed let $\phi\in C_0^\infty(B_\rho^+)$, $\rho$ is fixed and
$\rho<R$ then  $v_{L}^R$ is a minimizer and we have
$$J(v_{L}^R, B_\rho^+)\leq J(v_{L}^R+\phi, B_\rho^+),\qquad \forall  B_\rho^+\subset\R^n_+.$$
More explicitly it can be rewritten as
$\int_{B_\rho^+}\Lambda\X{v_L^R>0}\leq
\int_{B_\rho^+}2\na v_L^R\cdot \na \phi+|\na \phi|^2+\Lambda\X{v_L^R+\phi>0}$.

\smallskip

By a customary  compactness  argument and weak  convergence of  gradients we get
$$J(v_{L}, B_\rho^+)\leq J(v_{L}+\phi, B_\rho^+),\qquad \forall\phi\in
C_0^\infty(B_\rho^+).$$
The  same  argument leads  to  the existence of  $v_{S}-$the smallest global
homogeneous solution. Thus
$$v_{S}\leq u\leq v_{L},\qquad  \forall u\in\p_\infty.$$
\smallskip

Since the class $\mathcal P_\infty$ is scale and $e_3, \dots,e_n$ translation invariant
it follows that $v_S, v_L$ are homogeneous and depend only on $x_1$ and $x_2$ variables.

\medskip
Now let us explicitly compute $v_{L}$ and $v_{S}$. For this we write the Laplacian in
in polar coordinates
\setlength\arraycolsep{2pt}
\begin{eqnarray*}
\Delta w=\frac{1}{r}\left[\frac{\partial(rw_r)}{\partial r}+
\frac{\partial}{\partial\phi}\left(\frac{w_\phi}{r}\right)\right]=\frac{1}{r}\left[g(\phi)+g''(\phi)\right],
\end{eqnarray*}
where $w=rg(\phi)$. Recall that $v_S, v_L$ are harmonic  outside of the
zero set by Proposition \ref{prop}. This implies that $g$ is a linear combination
of $\sin \phi$ and $\cos \phi$.
Therefore the largest and smallest solutions are linear combinations of $x_1$ and $x_2$.

\smallskip
Assume that
\begin{eqnarray*}
  v^+=ax_1+bx_2,\ \  \textrm{in}\ \ \Omega^+(v), \qquad v^-=Ax_1+Bx_2,\ \ \textrm{in}\ \ \Omega^-(v),
\end{eqnarray*}
where $v^+$ and $v^-$ are respectively the positive and negative parts of $v$ and $a, b, A, B$ are constants to be determined.
The boundary condition $v=\alpha_+x_2^+-\alpha_-x_2^-$ on $\Pi$ implies
$b=\alpha_+, B=\alpha_-.$
\smallskip

Let us assume that the free boundary $\Gamma(v)$ is given by
$$x_1=x_2\tan\theta.$$
Both $v^+$ and $v^-$ must vanish on $\Gamma(v)$.
Hence
\begin{eqnarray*}
  0=ax_1+bx_2 = ax_2\tan\theta+\alpha_+x_2=x_2(a\tan\theta+\alpha_+)
\end{eqnarray*}
and we easily find that
$\displaystyle a=-\frac{\alpha_+}{\tan\theta}=-\alpha_+\cot\theta.$
Similarly
$$A=-\frac{\alpha_-}{\tan\theta}=-\alpha_-\cot\theta.$$
Summarizing we have
\begin{eqnarray*}
  v^+=\alpha_+(-{x_1}\cot\theta+x_2),\qquad  v^-=\alpha_- (-{x_1}\cot\theta+x_2).
\end{eqnarray*}
Note that $\cot\theta$ takes only two values, positive and negative, corresponding
respectively to large and small solutions.
To evaluate $\cot\theta$ we need to use the gradient jump condition
$|\nabla v^+|^2-|\nabla v^-|^2=\Lambda$,
which is now satisfied in classical sense, see Proposition \ref{prop}.
Substitution of $v$ into this identity gives
\begin{eqnarray*}
  \alpha_+^2(1+{\cot^2\theta})-\alpha_-^2(1+{\cot^2\theta})=\Lambda
\end{eqnarray*}
or equivalently
$${\cot\theta}=\pm\sqrt{\frac{\Lambda}{\alpha_+^2-\alpha_-^2}-1}.$$

Note that if $\Lambda\leq{\alpha_+^2-\alpha_-^2}$ then there is no free boundary.
Summarizing we get that
\begin{eqnarray}\label{globsol}
  v_{L}&=&\alpha_+(\gamma x_1+x_2)^+  -\alpha_-(\gamma x_1+x_2)^-,\\\nonumber
  v_{S}&=&\alpha_+(-\gamma x_1 +x_2)^+ -\alpha_-(-\gamma x_1+x_2)^-,\\\nonumber
  \gamma&=&\sqrt{\frac{\Lambda}{\alpha_+^2-\alpha_-^2}-1}.
\end{eqnarray}

The above discussion is  summarized in  the following proposition.

\begin{prop}\label{large-small-prop}
   The largest and smallest solutions $v_L, v_S$ are given by (\ref{globsol})
and these are the only two dimensional homogeneous global solutions.
\end{prop}

\subsection{Comparison of $W$-energy}
The aim of this section is to show that $v_S$ has bigger $W$-energy than $v_L$.
For all values of $\alpha_\pm$ for which $v_S \neq v_L$ we have

\begin{equation}\label{W-comp}
 W(1,v_{S}, 0) > W(1,v_{L}, 0),
\end{equation}

As a consequence we get that the largest solution is stable in the following sense:
\begin{prop}
 Let $u\in\mathcal P_1(n, \lambda_\pm, \alpha_\pm , g)$
and suppose there is $R_0\in (0, 1)$ such that $W(R_0, u, 0)< W(1, v_S, 0)$
then any blow-up limit $u_0$ of $u$ coincides with $v_L$.
\end{prop}

\pr To check this we recall the monotonicity of $W$,
to infer that $W(0^+, u, 0)=W(1, u_0, 0)< W(1, v_S, 0)$,
which in view of Theorem C
and $W(1,v_{S}, 0)\geq W(1,v_{L}, 0)$ implies that $u_0=v_L$.

Now it remains to show (\ref{W-comp}). If $v$ is a homogeneous solution,
then $W$ is constant hence it suffices to compute $W(1, \cdot, 0)$.
By Green's formula
$$\int_{B_1^+}|\nabla v|^2=\int_{\partial B_1^+}v\frac{\partial v}{\partial \nu}.$$
We can easily compute
\begin{eqnarray}\nonumber
   W(1,v)&=&\int_{\partial B_1^+}v\frac{\partial v}{\partial \nu}+\int_{B_1^+}\Lambda\chi\{v>0\}-\int_{S^+_1}v^2\\\nonumber
   &=&\int_{B_1'}v\frac{\partial v}{\partial
   \nu}+\int_{B_1^+}\Lambda\chi\{v>0\},
\end{eqnarray}
where the last equality follows from $v(x)=x\cdot \nabla v(x)$ on  $S_1^+=\partial
B_1^+\cap\R^n_+$.
In particular one can take $v$ to be $v_{L}$ or $v_{S}$.

Now let $\theta\in(0, \pi/2)$ be determined  from
$$\cot\theta=\sqrt{\frac{\Lambda}{\alpha_+^2-\alpha_-^2}-1}.$$
Utilizing the explicit form  of $v_S$ one can readily verify that
$$\int_{B_1^+}\Lambda\chi\{v_{S}>0\}=\int_0^1\int_{S_\rho^+}\Lambda\chi\{v_{S}>0\}=\frac1n\int_{S_1^+}\Lambda\chi\{v_{S}>0\}=
\frac{\Lambda\theta}{2 \pi}\omega_n, $$
where $\omega_n$ is the volume of $n$-dimensional unit ball.
Similarly
$$\int_{B_1^+}\Lambda\chi\{v_{L}>0\}=\frac{\Lambda(\pi-\theta)}{2 \pi}\omega_n.$$

\smallskip

Next  we notice that $\frac{\partial v}{\partial \nu}=-\frac{\partial v}{\partial x_1},$
on $B_1'$, therefore we have
\begin{eqnarray}\nonumber
  -\int_{B_1'}v_S\frac{\partial v_S}{\partial x_1}&=&-\int_{B_1'\cap\{x_2>0\}}\alpha_+x_2(-\gamma\alpha_+)-
   \int_{B_1'\cap\{x_2<0\}}\alpha_-x_2(-\gamma\alpha_-)\\\nonumber
   &=&\gamma\alpha_+^2\int_{B_1'\cap\{x_2>0\}}x_2+\gamma\alpha_-^2\int_{B_1'\cap\{x_2<0\}}x_2\\\nonumber
   &=&\gamma(\alpha_+^2-\alpha_-^2)\int_{B_1'\cap\{x_2>0\}}x_2\\\nonumber
   &=&\frac{\omega_{n-2}}{n}\gamma(\alpha_+^2-\alpha_-^2),
\end{eqnarray}
where $\gamma=\cot\theta=\sqrt{\frac{\Lambda}{\alpha_+^2-\alpha_-^2}-1}$.
Hence
\begin{eqnarray*}
  W(1,v_{S}, 0)=\gamma(\alpha_+^2-\alpha_-^2)\frac{\omega_{n-2}}{n}+\Lambda\frac{\theta\omega_n}{2\pi},
\end{eqnarray*}
and similarly one can see that
\begin{eqnarray}\nonumber
  W(1,v_{L}, 0)=-\gamma(\alpha_+^2-\alpha_-^2)\frac{\omega_{n-2}}{n}+\Lambda\frac{(\pi-\theta)\omega_n}{2\pi}.
\end{eqnarray}

Summarizing we have that
$$W(1,v_{S}, 0)- W(1,v_{L}, 0)=\Lambda\frac{\omega_{n-2}}{n}\left[2\gamma\frac{\alpha_+^2-\alpha_-^2}{\Lambda}+
\frac{(2\theta-\pi)}{2\pi}\frac{n\omega_n}{\omega_{n-2}}\right].$$
Using the explicit computation for $\omega_n$
we we obtain
\begin{eqnarray}\nonumber
 n\frac{\omega_n}{\omega_{n-2}}=2\pi.
\end{eqnarray}
Finally we observe that $\sin^2\theta=\displaystyle\frac{\alpha_+^2-\alpha_-^2}{\Lambda}$
 hence
\begin{eqnarray}\nonumber
 2\gamma\frac{\alpha_+^2-\alpha_-^2}{\Lambda}+
\frac{(2\theta-\pi)}{2\pi}\frac{n\omega_n}{\omega_{n-2}}&=&2\gamma\frac{\alpha_+^2-\alpha_-^2}{\Lambda}+
{(2\theta-\pi)}\\\nonumber
&=&2\left[\cot \theta\sin^2\theta-(\theta-\frac\pi2)\right]\\\nonumber
&=&\sin(2\theta)-2\theta+\pi\\\nonumber
&\geq 0&.
\end{eqnarray}

\noindent Therefore
$$W(1,v_{S}, 0)\geq W(1,v_{L}, 0)$$
and equality holds if and only if $\theta=\pi/2$.\qed


\medskip
\section{Proof of Theorem C }

\subsection{Free boundary as generalized minimal surface}
The aim of this section is to classify  homogeneous global solutions.
For $n=2$ this was done in Proposition \ref{large-small-prop}.
Therefore from now on we shall assume $n>2$, $\alpha_-=\lambda_-=0$ (i.e. the one phase case).
Notice that the condition $\lambda_-=0$ can be dropped due to the formulas (\ref{J-redef0})
 and (\ref{J-redef}).
We recall that if $u$ is a global solution, and hence local minimizer,
of $J$ for one phase problem then
\begin{equation}\label{grad-estim}
\sup_{B_r(x_0)}|\na u|^2\leq \Lambda +C(x_0)r^\alpha
\end{equation}
for any $x_0\in \R^n_+$, $B_r(x_0)\subset \R^n_+$ and $C(x_0)$ depends on
$\dist(x_0, \Pi)$, see Theorem 6.3 \cite{AC}.
As a result we obtain that for any free boundary
point $x_0$ the estimate holds
\begin{equation}\label{limsup-est}
 \limsup_{\begin{subarray}
            xx\in \Omega^+(u)\\ x\rightarrow x_0
           \end{subarray}
}|\na u(x)|^2\leq \Lambda .
\end{equation}
Our first task is to show that the estimate (\ref{limsup-est}) holds in ${\supp u}$.

\begin{lemma}\label{grad-estim-0}
 Let $u$ be a global homogeneous solution. Then for  $z^0 \in \Gamma  \cup \{x_1=0, x_2 >0\}$
$$
 \limsup_{\begin{subarray}
 {}z\in\R^n_+\\
 z\rightarrow z^0
          \end{subarray}
}\frac{|u(z) -u(z^0)|}{|z-z^0|}\leq \sqrt\Lambda.
$$
\end{lemma}

\proof
 To see this let $z_0\in \Pi$ and $u(z^0)>0$. Then there is $r>0$ such that  $u\in C^1(\overline{B_r^+(z^0)})$.
Thus the tangential derivatives are controlled by $\alpha_+\leq \Lambda$.
As for the normal derivative we notice that
from the definition of $v_S$ and $v_L$ we have that (since $\alpha_-=0$)
$$|\na v_S|=|\na v_L|=\sqrt\Lambda.$$
But $v_S\leq u\leq v_L$ and $v_S=u=v_L$ on $\Pi$, hence it is enough
to estimate the $x_1-$ derivative. Indeed, from the estimate
$v_S\leq u\leq v_L$ and $v_S(z^0)=v_L(z^0)=u(z^0)$ we get
$$\frac{\partial v_S(z^0)}{\partial x_1}\leq \frac{\partial u(z^0)}{\partial x_1}\leq
\frac{\partial v_L(z^0)}{\partial x_1}.$$
Therefore $|\na u(z^0)|^2\leq \Lambda$.

It is also apparent by the free boundary condition (\ref{grad-estim}) that $|\nabla u|^2 \leq \Lambda$ on the free boundary.
\qed

\begin{lemma}\label{limsup-est1}
 Let $u$ be a global homogeneous solution. Then
\begin{itemize}
\item[$\bf 1^\circ$] the following estimate is true
\begin{equation}
 \sup_{x\in \R^n_+ \cap \{u>0\}}|\na u(x)|^2\leq \Lambda.
\end{equation}
\item[$\bf 2^\circ$]
In particular $\Gamma(u)$ is a generalized surface of non-positive outward mean curvature.
\end{itemize}
\end{lemma}

It should be remarked that the estimate (\ref{limsup-est1})
does not hold for non-homogeneous global solutions; see \ref{non-convexity}.

\smallskip

{\bf Proof of Lemma \ref{limsup-est1}:}\
Suppose the statement of the lemma fails, then there is a maximizing sequence $x^j$ with the property that
$|\nabla u(x^j)|^2 \to \Lambda + \epsilon_0 > \Lambda $. By  zero-degree homogeneity of $|\nabla u |^2$ we may assume $x^j$
are on the unit sphere. Also by sub-harmonicity of $|\nabla u|^2$ we  assume that $x^j$ tend to the boundary of $\{u>0 \} \cap \{x_1>0\}.$
By Lemma \ref{grad-estim-0} the sequence $x^j$ cannot converge to either of the boundaries (free or fixed). Hence it converges to the "corner"-points
$\{x_1=x_2=0, \ |x|=1\}$.

Let $r_j=\dist(x^j, \Gamma \cup \Pi)$, then we have  three different possibilities:
\begin{eqnarray*}
 &&{\bf Case \ 1:}  \dist(x^j, \Gamma) \approx x^j_1 \Rightarrow \ r_j \approx  \dist(x^j, \Gamma) \approx x^j_1,\\
&&{\bf Case \ 2:} \dist(x^j, \Gamma) =o(\dist(x^j, \Pi)) \Rightarrow r_j = o(x_1^j),\\
&&{\bf Case \ 3}: \dist(x^j, \Pi)=o(\dist(x^j,\Gamma))\Rightarrow r_j= o(\dist(x^j, \Gamma)).
\end{eqnarray*}
Notice that $x^j_1=\dist(x^j, \Pi)$.
We shall see that all these cases will lead to a contradiction.

\smallskip

{\bf Case 1:} Let $\tilde x^j$ be the closest corner point on the $n-2$ dimensional unit sphere, i.e. $\tilde x^j \in \{x_1=x_2=0, \ |x|=1\}=\mathbb S^{n-2}$,
in first case, and in the other two cases the closest point  on the boundary to $x^j$ (we again assume this close point is on the unit sphere).

Now let $d_j=|x^j-\tilde x^j|$ and scale $u$ at   $\tilde x^j$   with $d_j$,
$$u_j(x)=\frac{u(\tilde x^j+d_jx)}{d_j} .$$
Note that $d_j\approx r_j\approx x_1^j$ translates to $u_j$ as follows; there is $y^j\in\mathbb S^{n}, y^j_3=\dots=y_n^j=0$
such that $y^j_1\approx \dist(y^j, \Gamma(u_j))\approx 1$ and
\begin{equation}\label{grad-e-0}
  \lim_{j\rightarrow \infty}|\na u_j(y^j)|^2=\Lambda+\e_0.
\end{equation}

Clearly  $u_j$ should be considered  in  a new domain, which is a scaled version of the support of $u$ at $\tilde x^j$
 and it contains $\supp v_S$.
In the two other cases below the support of $u_j$ converges to a half space.

Next we see that in all cases $u_j$ converges to a limit function $u_0$ (at least for a subsequence) with further property that $|\nabla u_0(y^0)|^2= \Lambda + \e_0$ (here $\displaystyle y^0=\lim_{j\rightarrow \infty} y^j$, again for a subsequence). In particular, and by construction, $ |\nabla u_0(x)|^2$ takes maximum at $y^0$, an interior point to the support of $u_0$. Hence by the strong maximum principle it must be constant, and therefore $|\nabla u_0(x)|^2= \Lambda + \e_0$ in the support of $u_0$. This in turn implies $u_0$ is linear.
But $u_0$ is a global minimizer, hence $|\na u_0|^2=\Lambda$ in $\supp u_0$ which in contradiction with (\ref{grad-e-0}).

\smallskip

{\bf Case 2:} Let
$u_j(x)=\frac{u(\tilde x^j+r_jx)}{r_j}$. We proceed as in Case 1 and extract a subsequence for which
$u_j\rightarrow u_0$ and $u_0$ is global minimizer. Furthermore (\ref{grad-e-0}) holds with $y^0=\lim y^j$ but in this
case $y^0\in \Gamma(u_0), y^0_3=\dots=y_n^0=0$. This implies that
$|\na u_0(y^0)|^2=\Lambda+\e_0$ which is in contradiction with (\ref{grad-estim-0}).

Now, in the first two cases, the free boundary is present (due to the length of scale $r_j$). In first case, we obtain a global minimizer in $\R^n_+$ with boundary data as before. At the same time we have $u_0$ is linear, which results into the fact that $u_0$ is one of the functions $v_L, v_S$. But then this contradicts the fact that
$|\nabla u_0 |^2= \Lambda + \e_0$.

\smallskip

{\bf Case 3:} Now the last case gives us scaling with center at the fixed boundary. Here we use both the small and the large solutions to bound the scaled function.
Indeed, for $\tilde x^j$ being the projection of $x^j$ onto $\Pi$, we have
$$
u_j(x)=\frac{u(r_jx+\tilde x^j) -u(\tilde x^j)}{r_j}= \frac{u(r_jx+\tilde x^j) - \alpha_+ \tilde x^j_2}{r_j}
$$
and hence the scaled versions of $v_S$ and $v_L$ at $\tilde x^j$ satisfy
$$
(v_S)_j \leq u_j \leq (v_L)_j .
$$
Hence the blow-up limits keep the order
\begin{equation}\label{sandwich}
v_S=(v_S)_0 \leq u_0 \leq (v_L)_0=v_L .
\end{equation}
Now as before we have $|\nabla u_0 |^2= \Lambda + \e_0$, and this is impossible due to (\ref{sandwich}), and the fact that $|\nabla v_L|^2=|\nabla v_S|^2= \Lambda $.

Now we turn to the proof of the second statement of Lemma \ref{limsup-est},
namely that $\Gamma(u_0)$ is a generalised surface of nonpositive outward mean curvature.
Let $S\subset \partial_{\rm{red}}\{u>0\}$ be a portion
of free boundary of $u$ and $S'$ a small perturbation of $S$ such that
$S'\subset \{u>0\}$ and $\partial S= \partial S'$. Then
$$\mathcal H^{n-1}(S)\leq \mathcal H^{n-1}(S')$$
i.e. $\partial_{\rm{red}}\{u>0\}$ is a generalized surface of non-positive outer mean curvature.
Notice that by Lemma \ref{lem-fin-perim} $\partial\{u>0\}$ has finite perimeter in
$B_1.$ Thus $\mathcal H^{n-1}(S)<\infty$.
\smallskip

To prove this we take the domains $G, G_0$ such that
$\partial G=S\cup S'$ and $\overline{G}\subset G_0\subset \R^n_+$.
Then we have
\begin{eqnarray}
 0=\int_G\Delta u=\int_S\partial_\nu u+\int_{S'}\partial_\nu u.
\end{eqnarray}
On $S$ we have that $\partial_\nu u(x)=|\na u(x)|=\sqrt{\Lambda}$,
for $\mathcal H^{n-1}$-a.e. $x\in\Gamma(u)\cap {G_0}$ \cite{AC}, whereas on
$S'$,  $|\na u|\leq \sqrt{\Lambda}$ by (\ref{limsup-est1}).
Comparing the integrals over $S$ and $S'$ we get that

$$\sqrt{\Lambda}\mathcal H^{n-1}(S)=\int_S\partial_\nu u=-\int_{S'}\partial_\nu u\leq \sqrt{\Lambda}
\mathcal H^{n-1}(S').$$
After canceling $\sqrt{\Lambda}$ the result follows.

\qed

\subsection{Preliminary Lemmas} Suppose that $u$ is a
third global homogeneous solution, which by Section
\ref{large-small-sect} satisfies $v_S\leq u\leq v_L$. In particular the
free boundary $\Gamma(u)$ lies  in between the planes $\Gamma_S$, and
$\Gamma_L$.

\smallskip

We first need a lemma that shows that free boundary is locally a graph.

\begin{prop}\label{normal-cont}
    Let $u$ be a global homogeneous minimizer and $\Gamma(u)$ touches tangentially
    the free boundary of $v_L$, at some point $x^0=(0,0,x_3^0, \cdots ,x_n^0) $ with $|x^0|=1$.
     Then in a small neighborhood of $x^0$ the free boundary $\Gamma(u)$
is a $C^1$ graph in the direction normal to $\Gamma_L$ in the upper half plane.
Moreover the normal vector to $\Gamma (u)$ is continuous up to the point $x^0$, and hence by homogeneity this holds on the axis $tx^0$ ($t>0$).
\end{prop}

Let $\Pi_0=\{x\in\R^n: x_1=x_2=0\}$ and $x^0\in\Pi_0\setminus\{0\}$
be any given free boundary point close enough to $\Pi_0$.
Let further  $\tilde x^0$ be the projection of $x^0$ onto $\Gamma_L$.
Then by tangential touch between the free boundary $\Gamma(u)$
and $\Gamma_L$ (which is a flat plane) one has that
$|x^0 - \tilde x^0|= o( x^0_1)$. In particular for $r_0=x_1^0$, sufficiently small,
we have that, in the ball $B_{r_0}(\tilde x_0)$,  the free boundary $\Gamma(u)$
is flat enough to satisfy the hypothesis of Theorem 8.1 in \cite{AC}.
In particular, in the direction of the plane $\Gamma_L$ the
free boundary is a $C^1$ graph locally in $B_{\frac{r_0}4}(\tilde x_0)$.
From here it follows that $\Gamma(u)$,
seen from the plane $\Gamma_L$, is $C^1$ graph over
$\Gamma_L \cap B_{\frac{r_0}{8}}(\tilde x_0)$.

It is now elementary to show that the normal of $\Gamma (u)$
is continuous up to the point $x^0$. Indeed, if this fails,
then there is a sequence $x^j$ on the free boundary
with normal $\nu^j$ staying uniformly away from the normal $\nu^L$  of
$\Gamma_L$, $|\nu^j - \nu^L|> \e_0 >0$.
Scaling $u$ at $x^j$ with $r_j=\dist(x^j,\Gamma_L)$ we
have a limit global minimizer in $\R^n$ (observe that this is due to tangential touch).
On the other hand the free boundary will then become a
plane, on one side of a scaled version of the plane $\Gamma_L$,
 but with the normal at the origin being $\nu^0$, with $|\nu^0 - \nu^L|> \e_0 >0$.
This is impossible.
\qed

\begin{lemma}\label{Tang-touch}
  Let $u\in\mathcal{HP_\infty},\Gamma(u)=\partial\{u>0\}$. If $u\geq 0, \alpha_-=0$ then $\Gamma(u)$
  does not touch $\Gamma_{L}$ tangentially.
\end{lemma}

\pr We argue towards a contradiction.
Let $x^0$ be a point  where the free boundaries touch each other.
We consider two possible locations: in first $x_1^0>0$ and then $x_1^0=0.$

\textbf{Case 1:}  Let us suppose that  $\Gamma(u)$ touches
$\Gamma_{L}$ at $x^0$ and $x^0_1>0$.
To conclude that this is a contradiction we use the free boundary condition and Hopf lemma. Notice that
in order to use
Hopf's lemma, we need (at least $C^{1,Dini}$)
regularity of $\Gamma(u)$ near $x^0$.

It follows from
the one side flatness, and classical regularity result  of Theorem 8.1. in
\cite{AC}. Then we can
apply Hopf's maximum principle to infer
$$\frac{\partial (u-v_{L})}{\partial \nu}(x^0)>0.$$
which is a contradiction in view of the tangential touch  condition.

\smallskip

\textbf{Case 2:} We first choose a new
coordinates system such that in new coordinates
$y=(y_1, y_2, \dots, y_n)$ we have $\Gamma_L=\{y\in\R^n: y_n=0, y_{n-1}>0\}$
and $\{u>0\}\subset \{y\in\R^n: y_n<0\}$.
Now let us assume that $\Gamma(u)$ touches the free
boundary of the larger solution $v_L$ at $y^0\not =0$ and $y^0_1=0$.
Then by Proposition \ref{normal-cont} the free boundary is locally a smooth graph,
seen from the plane $\Gamma_L$.
In particular near $y^0$, the free boundary
can be represented as $y_n=h(y')$, $y'=(y_1, y_2, \dots, y_{n-1})$ and  that
$h$ is a subsolution to the minimal surface equation in the weak sense.

Indeed,   let $ \R^{n-1}_+=\{y\in \R^n: y_n=0, y_{n-1}>0\}$ and $\widetilde B\subset \R^{n-1}_+$
be a ball touching
$\partial \R^{n-1}_+$ at $y^0$. Then by Lemma \ref{limsup-est1}, $\bf 2^\circ$ the
surface area functional will increase, if we replace $h$ by $h_\e=h+\e \phi$
for any $\phi\in C_0^\infty(\widetilde B), \phi\leq 0$ and $\e>0$ is small.
This comparison yields
$$\mathfrak Mh(y')=\div\left(\frac{Dh(y')}{\sqrt{1+|\na h(y')|^2}}\right)\geq 0
\qquad {\rm weakly\ in}\ \widetilde B.$$
Thus we have
\begin{equation*}
\left\{
\begin{array}{lll}
\mathfrak M h(y')\geq 0 \qquad \ \ &\mbox{in } \widetilde B,\\
h(y')\leq 0 & \mbox{in} \ \widetilde B,\\
h(y^0)=0 & y^0\in \partial \widetilde B.
\end{array}
\right.
\end{equation*}
By Hopf's  principle
$$\frac{\partial h(y^0)}{\partial y_{n-1}}>0$$
which is in contradiction with the tangential touch of $\Gamma(u)$ and $\R^{n-1}_+$.
\qed

\smallskip

From Lemma \ref{Tang-touch}
we know that $\Gamma(u)$ cannot touch $\Gamma_L$.
Using this observation we can construct yet another global minimizer
$\widetilde u$ such that it  is two dimensional
and distinct from $v_L$ and $v_S$. This, however, will contradict
Proposition \ref{large-small-prop}, and the proof of Theorem C will finish.

Thus to complete the proof of Theorem C we need to
construct  $\widetilde u$. This is done by the next lemma.

\begin{lemma}\label{tilt}
Let
$\mathbb V_\epsilon=\{x\in \R^n_+ : x_2<-(\gamma-\epsilon)x_1\}$ for small $\epsilon >0$.
If
$\mathbb{ V_\epsilon}\subset \{u=0\}$
then there is a two dimensional global solution $\widetilde u$ which is
distinct from $v_L$ and $v_S$.
\end{lemma}

\smallskip

\pr Suppose $\mathbb{V_\epsilon}\subset \{u=0\}$ for some $\epsilon >0$. Then
we can construct a global solution $\widetilde u$ such that $\widetilde u\geq u$,
$\widetilde u$ is two dimensional
and $\Gamma(\widetilde u)\subset \R^n_+\setminus \mathbb{V_\epsilon}$.

For $r>0$ fixed and $x\in B_r^+$, we put
$g_r(x)=\sup\{ u(x+\ell T), T\in \R, \ell\in \mathbb S^{n-2} \}$ where
$$\mathbb S^{n-2}=\{\ell=(0, 0, \ell_3, \dots, \ell_n),\ell_3^2+\ell_4^2+\dots\ell_n^2=1 \}.$$

Let $w\in \mathcal P_r(n, \lambda_\pm, \alpha_+, 0, g_r)$, i.e.
$w$ is a local minimizer of $J(\cdot, B_r^+)$ with $w=g_r$ on $\partial B_r^+$ see Remark \ref{rem-p-r}.
From Lemma \ref{rem-max-min} we infer that $\widetilde u_r=\sup w$ is a local minimizer
and $\widetilde u_r\geq w$ for any $w\in \mathcal P_r(n,  \lambda_\pm,\alpha_+,0, g_r)$.
In particular $\widetilde u_r\geq u$ in $B_r^+$.

Taking  $r_j\rightarrow \infty$, we have
from Proposition \ref{holder-0}, that  there is a
subsequence $r_{j_k}$ such that $\widetilde u_{r_{k_j}}\rightarrow \widetilde u_0$ locally in
$H^1$ and $C^0$ and $\widetilde u_0\in \mathcal {P_\infty}$.
Because $\mathcal P_r(n,\lambda_\pm,\alpha_+,0, g_r)$ is
translation invariant for each
$\ell\in \mathbb S^{n-1}$, it follows that $\widetilde u_0$ is
two dimensional solution. The condition $\mathbb V_\epsilon \subset \{u=0\}$
translates to
$\widetilde u_0$ and we get that
$\Gamma(\widetilde u_0)\subset \R^n_+\setminus \mathbb V_\epsilon.$ Furthermore
$$\widetilde u_0\geq u.$$
Since $v_L$ and $v_S$ are the only two dimensional
homogeneous global solutions, we conclude that $\epsilon=0$,
see Proposition \ref{large-small-prop}.
\qed

\begin{remark}
It is noteworthy that the classification of global homogeneous solutions for the two-phase case
 would have been available if one already knew that the free boundary is regular.
 Indeed, if we a priori know that  the free boundary is regular, then one can
apply maximum principle to $|\nabla u^+|^2$ in the set $\{ u >0\}$, and find out
 that the maximum must be on the boundary (either free or fixed).
  Actually, an argument similar to that of the proof of Lemma \ref{limsup-est1},
would then result in the fact that maximum is
 exactly on the boundary.

 Suppose now the maximum is on the free boundary.
  Then at such a maximum point $x^0$
(which is a maximum point for both $|\nabla u^+|^2$ due to Bernoulli
boundary condition $|\nabla u^+|^2 = \Lambda + |\nabla u^-|^2$) one gets that
$ \partial_\nu |\nabla u^+ (x^0)|^2 < 0$, where $\nu$ is the unit normal  pointing inwards support of $u^+$. From here along with a possible regularity of free boundary it follows that $2u^+_\nu u^+_{\nu\nu} (x^0) < 0$, which along with  $ u^+_\nu (x^0) > 0 $  gives that $u^+_{\nu\nu} (x^0) <0$.
  By representation of Laplacian on the free boundary we get $0=\Delta u^+ = \Delta_S u^+ + H u^+_\nu + u^+_{\nu \nu} $, and
 since $\Delta_S u^+ =0$,  $u^+_\nu (x^0) >0$, and $u^+_{\nu\nu} (x^0) <0$ we arrive at $H(x^0)>0$.
A similar argument applied to $u^-$ gives us
the converse $H(x^0)<0$, and we shall have a contradiction, unless $|\nabla u|$ is constant.

Next suppose the maximum for $|\nabla u^+|^2$  is on the fixed boundary $x^0 \in \{x_1=0\}$.
Then we have by a similar argument $u_1 (x^0) u_{11}(x^0) <0$.
Now with a representation of the Laplacian on $\{x_1=0\}$ along with
linearity of the boundary data we have $0= Hu_1 + u_{11}$.
Since the fixed boundary is a flat surface we have $H=0$, and hence
$u_{11}=0$ on the fixed boundary. This contradicts $u_1 (x^0) u_{11}(x^0) <0$.
\end{remark}


\medskip
\section{Proof of Theorem D}

Now we are ready to produce the proof of Theorem D,
exhibiting the non-tangential behavior of the free boundary.

\subsection*{Non-uniform approach}
Take $u\in \p_1(n,\lambda_\pm, \alpha_+, 0, g)$ and let $u_0$ be a blow-up
of $u$ at the origin. Then by Propositions \ref{technical} and \ref{W-prop} $u_0\in \mathcal{HP_\infty}$.
From Theorem C, $u_0$ is either $v_S$ or $v_L$. Suppose that $u_0=v_S$.
Let us consider the cone
$$K_\sigma^+: = \left\{x\in \R^n_+,x_2>0,
\frac{ x_2}{\gamma+\sigma}<x_1<\frac{ x_2}{\gamma-\sigma}\right\},$$
for small  $\sigma>0$ (cf. (\ref{cones})).
Then  we claim that for each $ \sigma>0$
there exist a $r_\sigma>0$
such that for any $r\in(0, r_\sigma),$ the following holds
\begin{equation}\label{thm-inc-D}
 \Gamma(u) \subset B_{r}^+\bigcap K_\sigma^+.
\end{equation}

This would suffice to conclude the tangential touch, since the modulus of
continuity can be constructed by inverting the relation
$\sigma \  \to \ r_\sigma$.

\smallskip
Suppose (\ref{thm-inc-D}) fails. Then there is  a sequence of
free boundary points $x^j\in\Gamma({u}),
|x^j|\rightarrow 0,  u\in \mathcal P_1(n, \lambda_\pm, \alpha_+,0, g)$
such that $x^j\not \in K_\sigma^+$ for some fixed $\sigma>0$.

Set $r_j=|x^j|$ and consider the  limit of the sequence $ u_j(x)=\frac{u(r_jx)}{r_j}$.
In view of Theorem A, $ u_j$'s are bounded and therefore by Proposition \ref{prop} and Theorem B
for a subsequence
$u_{j_m}\rightarrow u_0\in \mathcal{HP_\infty}$.
Moreover  the sequence of points $y^j=x^j/|x^j|\in \partial B_1^+$  is such that
$y^j \not \in K_\sigma^+$ and again by compactness this leads to the
existence of $y^0\in \partial B_1^+ \setminus  K_\sigma^+$ such that $ u_0(y^0)=0$.

From  monotonicity formula of Weiss, Proposition \ref{W-prop}
 one can also show that $u_0\in \mathcal{ HP_\infty}$  (see Section \ref{Weiss}) and hence we
can invoke Theorem C to conclude that $ u_0 $ is  $v_S$.
This contradicts the fact that $y^0\in\partial B_1^+\setminus  K_\sigma^+$,
and the proof of the first part is completed.
The case when $u_0=v_L$ is treated analogously.
\smallskip
\subsection*{Uniform approach}
To show the uniformity in the second statement of Theorem D,
we shall argue in the same way as above, but let $u$ change during the scaling.
In other words we define $v_j(x) =\frac{u_j(r_jx)}{r_j} $ with
$u_j \in  \mathcal P'_1(n, \lambda_\pm, \alpha_+,0, g)$, i.e.
$\lim\limits_{r\rightarrow 0}W(r,u_j, 0)=W(1,v_S,0)$.
As above the scaled functions will converge to a global solution $v_0$,
but $v_0$ is not necessarily homogeneous, and this is the only difference between the two cases.

Nevertheless, the assumption that
$\lim\limits_{r\rightarrow 0}W(r,u_j, 0)=W(0^+, u_j, 0)=W(1,v_S,0)$ for fixed $j$
implies that $W(tr_j,u_j, 0)=W(1, \frac{u_j(tr_jx)}{tr_j}, 0)=W(t,v_j, 0)\geq W(1, v_S, 0)$
by monotonicity of $W$ (see Proposition \ref{W-prop}) and after having sent $t$ to zero. This yields

$$W(t,v_0, 0)= \lim_{r_j\rightarrow 0}W(tr_j,u_j, 0)=\lim_{r_j\rightarrow 0}W(t,v_j, 0)\geq W(1, v_S, 0),
$$
where $v_0$ is the global limit of a subsequence of $v_j$.
The first inequality follows from strong convergence of $\na v_j$ in $L^2$,
since $\na v_j$ is a bounded sequence in $L^\infty$ and hence we can apply Theorem 1 from \cite{Z} and
Proposition \ref{technical}
to a suitable subsequence $\{r_j'\}\subset\{r_j\}$.

Next, the blow-down of $v_0$, at  infinity, i.e.
consider the scaling $v_0(rx)/r$ with $r\to \infty$
which results in a new homogeneous global solutions $v_{00}$. From  monotonicity formula, Proposition
\ref{W-prop}, we have
$$W(1,v_S,0)
\leq W(t,v_{0}, 0) \leq  W(\infty ,v_{0}, 0)= W( 1,v_{00}, 0).$$
Since $v_{00}$ is homogeneous, we can apply Theorem C to conclude $v_{00}$ is either
$ v_L$, or $v_S$. By the energy comparison (\ref{W-comp})  we should then have
$v_{00}=v_S$. Therefore
$W(t, v_0, 0)=W(1, v_S, 0)$ for any $t>0$ hence $v_0$ is homogeneous
by Proposition \ref{W-prop}.
Now, as in the previous case, contradiction comes from the fact that
$y^0\in\partial B_1^+\setminus  K_\sigma^+$.
\qed


\medskip

\section{Proof of Theorem E}
\subsection{Instability}
The problem studied in this paper is highly unstable in the
sense that changing the boundary data, no matter how small,
may result in a different behavior of
the touch between the free and the fixed boundary.
This behavior was already alluded in Theorem D, where we could not prove uniform behavior for solutions
that touch tangentially $\Gamma({v_L})$,
at the same time that the uniformity worked well for the class
${\mathcal P}_r'(n, \lambda_\pm, \alpha_\pm, g)$.

\smallskip

To illustrate this phenomenon, take $\alpha_-=0$ and
consider the largest  homogeneous  global solution $v_L$ as in Theorem C.
Consider now the minimization problem in the upper half
ball using the restriction of suitably scaled $v_L$ on the
boundary of $B_1^+$ as boundary data.
Now we know that the function itself is a minimizer.
Next let us decrease the data on the plane $\Pi$
to $(\alpha_+-\e ) x_2^+$.
A minimizer $u^\e$ of the functional with this boundary value on
$\Pi$ will exists, say take the smallest minimizer with boundary
values  $u^\e \leq v_L$ on $S_1^+$, so that  $u^\e\leq v_L$.
In particular this means that the free boundary for this minimizer will not touch the origin.
Indeed, if it touches the origin  then we can blow up $u^\e$ at the origin,
since by Theorem A $u^\e$ has linear growth at the origin,
and obtain a global minimizer $u^\e_0$, with data $(\alpha_+-\e)x_2$ on $\Pi$.

Now from the classification of the homogeneous global solutions, Theorem C,
we must have that $u^\e_0 < v_L$, and thus $u^\e_0 = v_S^\e$ in $B_1^+$.

This means that the free boundary cannot touch the origin, for any $\epsilon>0$.
In particular, by Theorem 5.1 in \cite{KKS}, we must have that it touches the fixed
boundary tangentially at some point $x^0$ with $x^0_2<0$.

\subsection{Non-homogeneous global solutions}
In this section we show the existence of a global solution  which is non-homogeneous.
We follow a perturbation method used in \cite{A-Sh}.

Let $\alpha_-=0, 0<\alpha_+<1, \Lambda=1$ and set $f^\e =(\alpha_+ - \e)x^+_2$.
Now consider a minimizer of
our functional in $B_1^+$,
with admissible functions having boundary data $f^\e$ on $\Pi$ and
$({\alpha_+-\e})({-\gamma x_1+x_2})^+$ on $S_1^+$ where
$\gamma=\sqrt{\frac1{\alpha_+^2}-1}$.

Let $\gamma_\e=\sqrt{\frac1{(\alpha_+-\e)^2}-1}$ then from Theorem C
$v_L^\e=(\alpha_+-\e)(-\gamma_\e x_2+x_1)^+$ is the largest
global homogeneous solution with boundary values $f^\e$ on $\Pi$.
Notice that $\gamma_\e>\gamma, f^\e\leq \alpha_+x_2^+$ implying that
$v_L\geq (\alpha_+-\e)(-\gamma x_2+x_1)^+$ on $S^+_1$.
Consider the class of local minimizers
$$\mathcal K_\e=\{u \in H^1(B_1^+), u=(\alpha_+-\e)(-\gamma x_1+x_2)^+\ \textrm{ on } \partial S_1^+,
u \textrm{ is a local minimizer of } J\}.$$
Then from the results of
Section \ref{build} $u^\e=\inf\limits_{\mathcal K_\e}u$ is a minimizer.
Furthermore $u^\e\leq \min(u^\e, v_L)\leq v_L$.

\begin{figure}
 \begin{center}
  \includegraphics[scale=0.9]{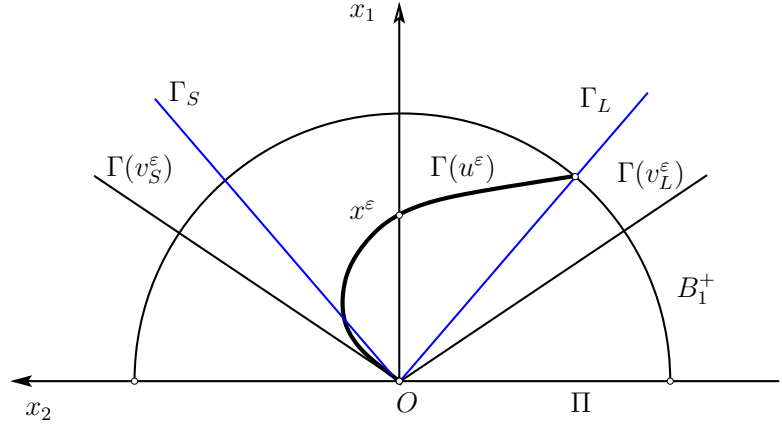}
\caption{The free boundary of $u^\e$}
 \end{center}
\end{figure}

For $\e$ fixed, any blow-up of $u^\e$ at origin is a
homogeneous global solution $u^\e_0 $,
which in view of the inequality  $u^\e \leq v_L$, implies $u_0^\e \leq v_L$.
Now $u^\e_0$ is a global homogeneous solution with boundary data $f^\e$,
and hence it must equal to one of  the functions
$(\alpha_+-\e)(\pm\gamma_\e x_2+x_1)^+$.
The only way for $u^\e_0$ to be as above and satisfy
$$
 u^\e_0 \leq v_L= \alpha_+(-\gamma x_2+x_1)^+
$$
is that $u^\e_0=(\alpha_+ -\e)(\gamma_\e x_2 - x_1)^+=v_S^\e$.
This in turn suggests that the free boundary $\Gamma({u^\e})$
starts at the origin with a tangential touch to $\Gamma({v^\e_S})$,
the smallest global solution with boundary data
$f^\e$.
Since the free boundary divides $B_1^+$ into two parts, it has to end
on $S_1^+$, see Figure 3.
In particular $\Gamma({u^\e})$
cuts the $x_1$-axis at some point $x^\e=(r_\e,0)$.
Now we consider the blow up of $u^\e$ with respect to $r_\e $.
Observe that $r_\e \to 0$
and thereby, utilizing Proposition \ref{technical} and choosing a suitable
subsequence, we obtain a global solution $u_0$ with
boundary data $\alpha_+ x_2$ and with $\Gamma(u_0)\ni (1,0)=\lim\limits_{\e\rightarrow 0} x^\e/r_\e $.
it follows from Theorem C that  this solution cannot be homogeneous.

\begin{remark}\label{non-convexity}
It should be remarked that in the above example of non-homogeneous global solutions, we have $|\nabla u|^2 \not \leq \Lambda $. Indeed, if this was true, then one may apply maximum principle to $|\nabla u|^2$ in $\{u>0\}$ and obtain a maximum on the free boundary (the free boundary is regular in $2$-space dimension).
Hence, by Hopf's lemma one obtains $\partial_\nu |\nabla u|^2 >0$, where $\nu $ is the unit  normal on the free boundary pointing outside the support of $u$.
In particular $u_\nu u_{\nu \nu} >0$. Since $u_\nu=|\nabla u|=\sqrt {\Lambda}$ we will have $u_{\nu \nu}>0$ on the free boundary. Using representation of Laplacian on
the free boundary $\Delta u = \Delta_S u + H u_\nu + u_{\nu \nu}$, where $H$ is the mean curvature, we conclude the convexity of the free boundary. This contradict the geometry of the example above.
\end{remark}

\section{Appendix 1}
In this section we prove that any blow up limit of
$u \in \mathcal P_r$ is homogeneous function of degree one.
The case when $g=0$ immediately follows from \cite{W}, Section 2.
When $g\not =0$ some extra care is needed, because the comparison of $u$
with its homogeneous extension $u_t(x)=\frac{|x|}tu(t\frac x{|x|})$ in $B_t^+$
fails on the flat portion of the boundary, i.e.
$u(x)\not =u_t(x)$ when $x\in \Pi\cap B_t$.

To fix the ideas we consider the model case $g(x)=C|x|^{1+\upkappa}$ with $\upkappa>0$
and $C=const$. Since $\rho^{-1}g(\rho x)\rightarrow 0$ as $\rho\downarrow 0$ it follows that
$u$ and $v=u-g$ have the same blowups at the origin.

\smallskip
\begin{lemma}\label{mon-aux}
 Let $u\in \mathcal P_r (n, \lambda_\pm, \alpha_\pm, g)$. Set $v=u-g$ where $g(x)=C|x|^{1+\upkappa}, \upkappa>0$. Then
$$\widetilde W(t)=\frac 1{t^n}\int_{B_t^+}|\na v|^2+ \Lambda\X{ v>-g}-
\frac 1{t^{n+1}}\int_{\partial B_t^+}v^2 +\frac{C_1}{\upkappa}t^\upkappa$$
is nondecresing function of $t$. Furthermore
\begin{eqnarray}
 \frac d{dt}\left\{\frac 1{t^n}\int_{B_t^+}|\na v|^2+ \Lambda\X{ v>-g}-
\frac 1{t^{n+1}}\int_{\partial B_t^+}v^2 +\frac{C_1}{\upkappa}t^\upkappa\right\}
\\\nonumber
\geq \frac 1{t^n}\int_{\partial B_t^+}\left(\na v\cdot \nu -\frac v t\right)^2.
\end{eqnarray}
\end{lemma}

\pr
Let $\phi\in H_0^1(B_r^+), r\in(0, 1)$ and let us define
$v=u-g, v_\phi=u+\phi-g$ in $B_1^+$. Then  $J(u)\leq J(u+\phi)$ transforms into
$J(v+g)\leq J(v_\phi+g)$. Employing Green's identity
we obtain
\begin{eqnarray*}
 \int_{B_r^+}|\nabla u|^2&=&\int_{B_r^+}|\nabla v|^2+2\int_{B_r^+}\nabla v \cdot \nabla g
+\int_{B_r^+}|\nabla g|^2\\\nonumber
&=&\int_{B_r^+}|\nabla v|^2-v(2\Delta g)+2\int_{\partial B_r^+}v  (\nabla g\cdot \nu)
+\int_{B_r^+}|\nabla g|^2.
\end{eqnarray*}

Utilizing this computation and the fact $v_\phi-v=H^1_0(B_r^+)$ we see
that if $u$ is a minimizer of $J(u, B_r^+)$, subject to its own
boundary values on $\partial B_r^+$, then $v$ is a minimizer of
\begin{equation}
 \widetilde J(v)=\int_{B_r^+} |\na v|^2-v(2\Delta g)+\Lambda\X{v>-g},
\end{equation}
because $2\int_{\partial B_r^+}w  (\nabla g\cdot \nu)
+\int_{B_r^+}|\nabla g|^2$ is constant for any
$w\in H^{1}(B_r^+), w|_{\partial B_r^+}=v|_{\partial B_r^+}$.

Thus it remains to prove that any blow up limit of $v$ at the origin
is a homogeneous function of degree one.

Let $t>0$ be small and take $\displaystyle v_t(x)=\frac{|x|} tv(t\frac x{|x|})$,
then on $\partial B_t$ $v_t$ agrees with $v$ and it follows $\widetilde J(v)\leq \widetilde J(v_t)$.
Using the homogeneity of $v_t$ and the identities
\begin{eqnarray}\nonumber
\nabla v_t(x)&=&\frac x{t|x|}v(t\frac x{|x|})+\nabla v(t\frac x{|x|})-\na v(t\frac x{|x|})\cdot \frac
x{|x|}\frac x{|x|}, \\\nonumber
|\na v_t|^2&=&t^{-2}v^2(t\frac x{|x|})+\left|\na v(t\frac x{|x|})\right|^2-
\left(\na v(t\frac x{|x|})\cdot \frac
x{|x|}\right)^2
\end{eqnarray}
one can easily compute

\begin{eqnarray}\nonumber
\int_{B_t^+}|\na v_t|^2+ \Lambda\X {v_t>-g}
&=&
\int_{B_t^+}\left[\frac x{t|x|}v(r\frac x{|x|})+\na v(t\frac x{|x|})-\na v(t\frac x{|x|})\cdot \frac
x{|x|}\frac x{|x|}\right]^2+\Lambda \X {v_t>-g}\\\nonumber
&= & \int_0^t \int_{\partial B_\rho^+}
\left[ t^{-2}v^2(t\frac x{|x|})+\left|\na v(t\frac x{|x|})\right|^2-
\left(\na v(t\frac x{|x|})\cdot \frac
x{|x|}\right)^2\right]+\Lambda\X{v_t>-g}\\\nonumber
&=&\frac  tn \int_{\partial B_t^+}|\na v|^2  +\frac  tn \int_{\partial B_t^+}
\left[\frac{ v^2}{t^2} -(\na v\cdot \nu)^2\right]+\int_{B_t^+}\Lambda\X{v_t>-g}.\\\nonumber
\end{eqnarray}

To deal with the last integral, we first notice that
$\{v_t(x)>-g(x)\}\subset \{v(t\frac x{|x|})>-Ct^{1+\upkappa}\}$.
Indeed if $x\in \{v_t(x)>-g(x)\} $ then $\frac{|x|}tv(t\frac{x}{|x|})>-C|x|^{1+\upkappa}$,
or equivalently $v(t\frac{x}{|x|})>-Ct|x|^\upkappa$. But $|x|\leq t$ since
$x\in B_t^+$. Thus $-Ct|x|^{\upkappa}\geq -Ct^{1+\upkappa}=-g(t)$.
In particular we get that
$\int_{B_t^+}\Lambda\X{v_t>-g}\leq \int_{B_t^+}\Lambda\X{v(t\frac x{|x|})>-g(t)}$
which, after applying Fubini's theorem, yields
$$\int_{B_t^+}\Lambda\X{v(t\frac x{|x|})>-g(t)}=\frac tn\int_{\partial B_t^+}\Lambda\X{v>-g}.$$

Next we notice that if $w\in H^{1}(B_t^+)$ and $|w(x)|\leq C|x|, x\in B_t^+$
then $\left|\int_{B_t^+}w(2\Delta g)\right|\leq C_1t^{n+\upkappa}$ with some tame constant $C_1$.
Therefore comparing the $\widetilde J$ energies in $B_t^+$, we get

\begin{eqnarray}
 0 &\leq & \widetilde J (v_t)- \widetilde J (v)\\\nonumber
&\leq &  \frac  tn \int_{\partial B_t^+}|\na v|^2 +\Lambda\X{v>-g} +\frac  tn \int_{\partial B_t^+}
\left[\frac{ v^2}{t^2} -(\na v\cdot \nu)^2\right] \\\nonumber
&&-\int_{B_t^+}|\na v|^2+ \Lambda\X{ v>-g}+C_1t^{n+\upkappa}\\\nonumber
&\leq &  \frac  tn \int_{\partial B_t^+}|\na v|^2 +\Lambda\X{v>-g}-
\int_{B_t^+}|\na v|^2+ \Lambda\X{ v>-g}\\\nonumber
&& +\frac  tn \int_{\partial B_t^+}
\left[\frac{ v^2}{t^2} -(\na v\cdot \nu)^2\right]+C_1t^{n+\upkappa}\\\nonumber
&=&\frac {t^{n+1}} n\frac d{dt}\left\{\frac 1{t^n}\int_{B_t^+}|\na v|^2+ \Lambda\X{ v>-g}\right\}\\\nonumber
&&-\frac t n\int_{\partial B_t^+}\left(\na v\cdot \nu -\frac v t\right)^2 - \frac  {2t}n \int_{\partial B_t^+}
\frac{v}{t^2}\left[\na v\cdot \nu-\frac{v}{t}\right]+C_1t^{n+\upkappa}.
\end{eqnarray}
Multiplying both sides by $nt^{-n-1}$  we conclude

\begin{eqnarray*}
 \frac d{dt}\left\{\frac 1{t^n}\int_{B_t^+}|\na v|^2+ \Lambda\X{ v>-g}-
\frac 1{t^{n+1}}\int_{\partial B_t^+}v^2 +\frac{C_1}{\upkappa}t^\upkappa\right\}
\\\nonumber
\geq \frac 1{t^n}\int_{\partial B_t^+}\left(\na v\cdot \nu -\frac v t\right)^2.
\end{eqnarray*}
\qed

\begin{remark}
{This argument shows that $g$ can be replaced by any
homogeneous polynomial or function of degree $m>1$ }.
\end{remark}

\begin{corollary}
 Let  $v$ be as in Lemma \ref{mon-aux}. Then any blow up limit of $v$ at the origin is
homogeneous of defree one. In particular any blow up of $u$ is homogeneous of degree one.
\end{corollary}
\pr The first statement follows exactly as in \cite{W}, Section 2. To show that the blow up of $u$
is homogeneous we need to notice that $g(rx)r^{-1}\rightarrow 0$ uniformly as $r\rightarrow 0$.
Hence the blow ups of $u$ and $v$ coincide.

\section{Appendix 2}

We shall discuss the rectifiablity of the free boundary in $B_1$.

\begin{lemma}\label{lem-nondeg}
 Let $u$ be a global homogeneous minimizer and $\Gamma(u)$ touches tangentially
the free boundary of $v_L$, then $u$ is nondegenerate, i.e. there is a
tame constant $c>0$ such that for any $x\in \Gamma(u)$ the following estimate is true
\begin{equation}
 \sup_{B_r^+(x)}u \geq cr, \qquad \forall B_r(x)\subset \R^n.
\end{equation}

\smallskip

\begin{remark}\label{rem-int-vs-sup}
 In \cite{AC}, a different form of nondegeneracy is proven
(see Lemma 3.4 in \cite{AC}), namely
 $$\fint_{\partial B_r(x)}u\geq cr, \qquad x\in \partial \{u>0\}.$$
This integral inequality implies that there is $y\in \partial B_r(x)$ such that
$u(y)\geq cr$. But $u$ is subharmonic, therefore $\sup_{B_r(x)}u\geq cr$.
\end{remark}

\end{lemma}
\pr Let $x\in \partial\{u>0\}$ and set $\delta(x)=\dist(x, \Pi)$. If
$\delta(x)\geq r$ then $B_r(x)\subset \R^n_+$. Taking
$u_r(y)=\frac{u(x+ry)}{r}, y\in B_1$ and employing Lemma \ref{limsup-est1} $\bf 2^\circ$
we see that $u_r$ is a local minimizer. Hence from remark \ref{rem-int-vs-sup} we obtain
$\sup_{B_{\frac 12}}u_r\geq c$, which after scaling back implies the desired result.

\smallskip

\noindent Now assume that $\delta(x)<r$. We consider two possible scenarios:

{\bf Case a)} $\frac{r}{1000}\leq \delta(x).$ Then using Remark \ref{rem-int-vs-sup}
in $B_{\delta(x)}(x)$ we get
$$\sup_{B_r^+(x)}u\geq \sup_{B_{\delta(x)}}u\geq \frac{c\delta(x)}{2}
\geq \frac{cr}{2000}.$$
\smallskip

{\bf Case b)} $\delta(x)<\frac{r}{1000}$. Let $R(x)=\dist(x, \Pi_0)$, where
$\Pi_0=\{x\in\R^n: x_1=x_2=0. \}$ and take $x_0\in\Pi_0$ such that
$R(x)=|x-x_0|.$ Notice that $R(x)\sim \delta(x)$, because $\Gamma(u)$ touches
$\Gamma_L$ tangentially. This means that there are two positive
constants $a, b$ such that $aR(x)\leq \delta(x)\leq bR(x)$  if
$x$ is close to $\Pi$ (see definitions of the cones $K_\sigma$).
We have $r>1000\delta(x)\geq a1000R(x)$ yielding $R(x)\leq \frac{r}{a1000}$.
In particular $$\rho=r-R(x)\geq r-\frac{r}{a1000}\geq \frac{r}{100}.$$

Observing that $B_\rho^+(x_0)\subset B_r^+(x)$ we get
$$\sup_{B_r^+(x)}u\geq
\sup_{B_\rho^+(x_0)}u\geq \sup_{B_\rho^+(x_0)}u_S\geq c\rho\geq \frac{cr}{1000}.$$\qed

\smallskip

\begin{lemma}\label{lem-fin-perim}
 Let $u$ be as in Lemma \ref{lem-nondeg}. Then
$$\mathcal H^{n-1}(B_1\cap \partial\{ u>0\})<\infty.$$
\end{lemma}

\pr For each open ball $B_r(x)\subset \R^n$ let $B_r^+(x)=B_r(x)\cap \R^n_+$.
Introduce the measure $\mu=\Delta u$. Clearly $\mu$ is nonnegative Radon measure, because
$\int_{B_r^+(x)}\mu=\int_{\partial B_r^+(x)}\na u\cdot \nu\leq Cr^{n-1}$.
Hence for any compact $D\subset \R^n$ we can cover $\overline{D\cap\R^n_+}$
by a finite number of balls, which yields $\mu(\overline{D\cap\R^n_+})<\infty$.

\smallskip

Next we want to show that there is a positive
constant $c_0$  such that for each $x\in \overline{B_{1}^+\cap \Gamma(u)}$
we have
\begin{equation}\label{lap-nondeg}
 \int_{B_r^+(x)}\mu\geq c_0r^{n-1}
\qquad {\rm if}\ r>0\ {\rm is\ small}.
\end{equation}
From (\ref{lap-nondeg}) one can conclude the proof
of Lemma by employing a standard covering argument.

First we note that by Lemma \ref{lem-nondeg} $u$ is nondegenerate,
that is there is a constant
$c>0$ such that
\begin{equation}\label{nondeg-app}
 \sup_{B_r^+(x)}u\geq cr,\qquad  \forall x\in B_{1}^+\cap \Gamma(u)
\end{equation}
for small $r>0$.

Now suppose that (\ref{lap-nondeg}) fails. Then
there is a sequence of free boundary points $x_j\in \Gamma(u)$
and a sequence of positive numbers $r_j>0$ such that
\begin{equation}\label{lap-fail}
 \fint_{B_{r_j}(x_j)}\mu\leq \frac{r_j^{n-1}}{j}.
\end{equation}

First, let us suppose that there is a subsequence $r_{j(m)}$ such that
$B_{r_j(m)}\cap \Pi_0\not =\emptyset$. Let $x_j^0\in \Pi_0$
and $\dist(x_j, \Pi_0)=|x_j-x_j^0|$.
Then consider $v_m(x)=\frac{u(x_j^0+r_{j(m)}x)}{r_{j(m)}}, x\in B_2^+$.
From Proposition \ref{technical} and Lemma \ref{limsup-est1} we get
$v_{m_k}\rightarrow v_0, \mu_{m_k}\rightharpoonup \mu_0$ where
$\mu_{m_k}=\Delta v_{m_k}$ and $\Delta v_0=\mu_0$, at least for a subsequence $m_k$, and
$v_0$ is a local minimizer. Moreover (\ref{lap-fail}) translates to
$$\fint_{B_1^+(y^0)}\mu_0=0$$
for some $y^0\in \Gamma_L$, i.e. $v_0$ is harmonic in
$B_1^+(y^0)$. From the strong maximum principle we conclude
$v_0=0$ which is in contradiction with nongedeneracy of $v_{m_k}$ and $v_0$.

\smallskip

Finally let us assume that $B_{r_j}(x_j)\cap \Pi_0=\emptyset$
for any $j$. Denote $\delta_j=\dist{(x_j, \Pi)}$.
From tangential touch of $\Gamma(u)$ and $\Gamma_L$ it follows that
$aR_j\leq \delta_j\leq b R_j$, where $R_j=\dist(x_j, \Pi_0)$. Thus we have
$r_j<R_j$. If, moreover, $r_j\geq \delta_j$ then
applying Theorem 4.3 \cite{AC} to $\frac{u(x_j+r_jx)}{r_j}$
we will conclude a contradiction if $j$ is large enough.

Thus without loss of generality we may assume that $\delta_j<r_j<R_j$.
Introduce $w_j(y)=\frac{u(x_j+\delta_j)}{\delta_j}, y\in B_1$ then
$$\fint_{B_{\delta_j}(x_j)}\mu \leq \fint_{B_{r_j^+}(x_j)}\mu\leq
\frac{r_j^{n-1}}{j}\leq \frac{\delta^{n-1}_j}{ja^{n-1}}.$$
Hence for $\Delta w_j=\mu_j$ we have $\int_{B_1}\mu_j\leq \frac{1}{ja^{n-1}}$.
On the other hand $\sup_{B_{\frac12}}w_j\geq c$. Extracting a subsequence
for which $w_j\rightarrow w_0, \Delta w_j\rightharpoonup \Delta w_0$ in $B_1$
at least for a subsequence,
where $w_0$ is a local minimizer in $B_1$, see Proposition \ref{technical}.
But $\fint_{B_1}\Delta w_j\leq \frac{1}{ja^{n-1}}\rightarrow 0.$
Thus $w_0\geq 0$ is harmonic and nondegenerate in $B_1$ and $w_0(0)=0$.
Hence by strong maximum principle $w_0=0$ which is in contradiction with
$\sup_{B_{\frac12}}w_0\geq c$.
\qed



\end{document}